\documentclass[12pt]{article}
\usepackage{amsmath,amssymb,array,graphicx,amsfonts,amsthm,subfig}
\usepackage{epsfig,algorithmic,algorithm,epstopdf,url}
\input epsf
\numberwithin{equation}{section}
\numberwithin{figure}{section}

\newtheorem{theorem}{Theorem}[section]

\def\OP{\Omega^+}
\def\OM{\Omega^-}
\def\d{\partial}

\newcommand{\ds}{\displaystyle}

\def\<{\langle}
\def\>{\rangle}

\newcommand{\X}{\ensuremath{\mathcal{X}}}

\author{Cary Humber\thanks{Naval Surface Warfare Center, Panama City Beach, FL and Department of Mathematics, North Carolina State University, Raleigh, NC} \,\, , 
Kazufumi Ito \thanks{ Center for Research in Scientific Computation, Department of Mathematics, North Carolina State University, Raleigh, NC} \,\, , Chad Bouton\thanks
{Battelle Memorial Institute, Columbus, OH}}

\title{Nonsmooth Formulation of the Support Vector Machine for a Neural Decoding Problem}

\date{}
\begin{document}
\maketitle

\begin{abstract}
This paper formulates a generalized classification algorithm with an application to 
classifying (or `decoding') neural activity in the brain.  Medical doctors and researchers have long been 
interested in how brain activity correlates to body movement.  Experiments have been 
conducted on patients whom are unable to move, in order to gain insight as to 
how thinking about movements might generate discernable neural activity.  Researchers 
are tasked with determining which neurons are responsible for different imagined movements 
and how the firing behavior changes, given neural firing data. For instance, 
imagined movements may include wrist flexion, elbow extension, or closing the hand. This is 
just one of many applications to data classification.  Though this article deals with an 
application 
in neuroscience, the generalized algorithm proposed in this article has applications in 
scientific areas ranging from neuroscience to acoustic and medical imaging.

\end{abstract}

\section{Introduction}
This paper deals with the problem of classifying neural activity
that correlates to specific imagined movements given data
recorded from an electrode array implanted in the primary motor cortex
of the human brain.
We are motivated by the training methods and neural decoding algorithms that have been developed by the Battelle Memorial Institute for 
the Braingate project as in \cite{bouton}.
The goal of the technology is to isolate and predict arm/hand movements a given patient is thinking about from signals processed by an
 electrode array. 
For this particular problem, we are trying to identify the active neurons, yet the data largely consists of measurements corresponding to the neurons at 
``rest".  In order to identify 
the 
neurons, we must first classify the data by movement type.  
We have 96 neural channels due to the construction of the Utah electrode array being utilized, and the training data contains neural firing reactions 
to a person simply thinking 
about the movements.
The overall goal of this neural network study is to identify the active neurons that are most responsible for each movement.  
Our objective is to have a simple linear classifier that incorporates sparsity in the formulation.  
Sparsity optimization is becoming increasingly utilized due to storage and implementation considerations (see for instance \cite{wright2,wright1}).
Not only that, but simple solutions are more robust to noise than exact solutions.
By sparsity, we mean 
a solution with the fewest number of nonzeros necessary to capture the core properties of the solution.
Mathematically, sparsity is represented
by $\ell_p$ norms for $0\le p\le 1.$
To this end, we obtain a simple linear classifier that, though simple, has desirable characteristics 
for classifying data.  In general, classification is a highly useful but difficult problem.  For certain applications, such
as the one considered in this paper, the difficulties can be compounded when the data is dynamic in time.  We address some necessary improvements
to widely used classifiers that only minimally increase computational effort.

Classifying the movements from the neural data can be formulated as a nonsmooth constrained minimization problem which leads to a, possibly 
ill-posed, inverse problem.  
There have been several methods proposed for classification problems of this type.  In this paper, we propose and implement several 
improvements to such classification 
methods.  
We propose a new formulation using sparsity for the weights of the classifier which yields a sharper classification
due to the noble use of the nonsmooth performance index (see \cite{regparam,newchoice}).
  Our formulation is based on the method proposed in \cite{cortes} and the improved least squares 
formulation developed in \cite{psvm}.  
Utilizing the $\ell_p$-norm for $0<p\le 1$ enables us to reduce the number of classified data and, hence,
reduce the misclassification of erroneous data.
In terms of the application considered in this paper, this approach improves the sparsity of the solutions 
in terms of the identified neurons for each movement. 
For this application, the Proximal Support Vector Machine(PSVM) was previously used due to the fact that it offers a 
closed form solution and has a computationally efficient implementation.
It will be demonstrated that our classification method attains solution sparsity, improves separability of the classes and is able to account 
for bias in the data. 
In order to improve the separability of the classes, we employ the 
penalty formulation of the inequality constraint in the linear Support Vector Machine(SVM). 
Due to the large amount of data for periods of rest, the standard classification algorithm could have a bias towards identifying the
 ``rest'' state, hence reducing
the identification of neurons corresponding to an imagined movement. To counteract this 
effect, we weight the error function accordingly.  
Moreover, we propose a new measure of performance for identifying the responsible neurons 
once the approximate solution 
has been computed.  The power of the methods presented in this paper is due to the fact
that not only is the data classified, but the strength of the classification in terms of the separation is obtained
from the same method.  That is, without extra work one obtains a measure of how the classifier performs along with the classifier.
We are motivated by the particular Braingate application, however the formulation developed here has the
potential for improving classification for a wide range of applications.

The paper is organized as follows.  In Section \ref{isvm}, we briefly discuss the SVM and PSVM, pointing out in what way
the methods can be improved.  Section \ref{newapproach} contains our approach for addressing the necessary improvements, which is 
largely based on sparsity optimization.  In Section 
\ref{choicerules} we present a choice rule for selecting the regularization parameter, which can be used for any of the methods discussed
in this paper.  Finally, in Section \ref{numres}, we present results from applying the PSVM and the sparse approach to real world
neural firing rate data.  

\section{The Support Vector Machine}
\label{isvm}


In this section, we give a basic outline of the Support Vector Machine (SVM) algorithms.
We are given training data $\mathcal{D}$, a set of $n$ points of the form
$$
\mathcal{D} = \left\{ ({x}_i, d_i)\mid {x}_i \in \mathbb{R}^m,\, d_i \in \{-1,1\}\right\}_{i=1}^n
$$
where the $d_i$ is either 1 or -1. We want to find the maximum-margin hyperplane
that divides the points having $d_i = 1$ from those having $d_i =-1$. Any hyperplane can be written as
the set of points $x$ satisfying
$$
x_i \cdot w - \gamma=0.
$$
To this end the linear SVM determines the hyperplane $(w,\gamma)^\intercal$ by
the constrained minimization;
\begin{equation} \label{Org}
\begin{array}{l}
\min \limits_{(w,\gamma)} \nu\,\sum\limits_{i=1}^m y_i+\dfrac{1}{2}(|w|^2+\gamma^2)
\\ \\
\mbox{subject to } d_i(x_i \cdot w-\gamma)\ge 1-y_i,\quad y_i \ge 0
\end{array}
\end{equation}
where $y_i$ measures the degree of misclassification and $\nu>0$
is a chosen parameter. That is, the SVM algorithm classifies data
into two categories, $\OM$ and $\OP$, geometrically separated by
the plane $\{x:x\cdot w=\gamma\}$, and clustered around the two planes
\begin{equation} \label{class}
\begin{array}{l}
\OM =  \{x \in \mathbb{R}^m :\; x\cdot w- \gamma \le -1\}
\\ \\
\OP =  \{x\in \mathbb{R}^m: \; x \cdot w-\gamma \ge +1 \}.
\end{array} \end{equation}
The classes and data $\mathcal{D}$  for this particular application will be discussed in Section
\ref{numres} to follow.

The authors of \cite{psvm} formulate the inequality $y \ge 0$ in terms of a penalty,
\begin{equation} \label{Iterm}
\min\limits_{w,\gamma,y}\quad \frac{\nu}{2}|y|^2+\dfrac{1}{2}(|w|^2+\gamma^2)\quad
\mbox{subject to  } D(Aw-\gamma e)+y \ge e,
\end{equation}
where $D=\operatorname{diag}(d_i)$ and $A \in \mathbb{R}^{n\times m}$ with rows
$A_i=x_i,\;1\le i\le n$. Furthermore, the PSVM algorithm in \cite{psvm} replaces the
inequality as the equality constraint and formulate the
unconstrained minimization
\begin{equation} \label{PSVM}
\min\quad \frac{\nu}{2}|y|^2+\frac{1}{2}(|w|^2+\gamma^2)\quad
\mbox{subject to  } D(Aw-\gamma e)+y = e.
\end{equation}
Our formulation is also an unconstrained minimization of the form
\begin{equation} \label{Our}
\min\quad \frac{\nu}{2}|y^+|^2+\frac{1}{2}(|w|^2+\gamma^2)\quad
\mbox{subject to  } D(Aw-\gamma e)+y = e,
\end{equation}
where $y_i^+=\max(0,y_i)$.
Our motivation for choosing $y^+$ in this manner can be understood by the following simple argument.
Note that if $y_i\ge0$ in \eqref{PSVM} then
$$
d_i(x_i\cdot w-\gamma)= 1-y_i
$$
and thus $y_i$ is the degree of misclassification.
However, if $y_i \le 0$
$$
d_i(x_i\cdot w-\gamma)=1-y_i \ge 1
$$
and thus the case is allowed. We are motivated by this fact to only penalize
$y^+_i=\max(0,y_i)$ in the formulation \eqref{Our}. In this
sense the formulation \eqref{Our} is penalizing the inequality
constraint of \eqref{Org}. It should improve the separability of the
classes based on the least squares formulation \eqref{PSVM}.

But, the advantage of \eqref{PSVM} is that it has the  closed form
solution $(w\;\gamma)^\intercal$. That is, \eqref{PSVM} is equivalent to
$$
\frac{\nu}{2}|H(w\; \gamma)^\intercal-e|^2+\frac{1}{2}(|w|^2+\gamma^2)
$$
where $H=D[A\;-e]$ (i.e. $y=e-H(w\; \gamma)^\intercal)$ and thus
$$
u=(w \;\gamma)^\intercal=(I+\nu\,H^\intercal H)^{-1}H^\intercal e.
$$
However, it will be shown in Section 
\ref{newapproach} that the formulation \eqref{Our} has an efficient implementation as well.

An important consideration for classification problems of this form is the possibility of ill-posedness.
If $H$ is very ill-conditioned, i.e. the singular values of $H$ decrease very rapidly
to zero, then the solution is very sensitive to the selection of $\nu>0$.
Thus, we need to develop a selection rule for $\nu$, which will be addresed in Section \ref{choicerules}
in a more general set-up.  The second term $\frac{1}{2}(|w|^2+\gamma^2)$ in \eqref{PSVM} represents
the $2$-norm of $u=(w\;\gamma)^\intercal$. It is more reasonable to use some other norms to
obtain a desirable classifier. One of our requirements is that fewer nonzero components of $w$ are in the final solution.
To satisfy this requirement, we use the $\ell^p$ norm with $0<p\le 1$
for our formulation in Section \ref{newapproach} to obtain the sparse solution.
The nonzero components of $w$ represent the essential and critical neurons for classifying the
specified movement.  In this way we can obtain the neural network information of the Braingate
technology. This point will be examined in the numerical tests presented in Section 
\ref{numres}.

\section{New Approach}
\label{newapproach}


Our general classifier can be written as the optimization problem
\begin{equation}\label{Gen}
\begin{array}{l}
\min\limits_{(w,\gamma)}\;\;\phi(y)+\beta\psi(w,\gamma),\\
\\
y=e-Hu
\end{array}
\end{equation}
where we assume $\phi,\psi$ are lower semi-continuous, as in
\cite{regparam},\cite{newchoice}. Typically $\phi$ and $\psi$
are chosen to be some norm on $\mathbb{R}^m$ (in general, a Banach space
$X$). For instance, in the PSVM we have $\phi(y)=|y|_2^2$ and
$\psi(w,\gamma)=\dfrac{1}{2}(| w|_2^2 +|\gamma |^2)$, for
$X=\mathbb{R}^n$, with $\beta=\frac{1}{\nu}$. As stated in the
previous section, the choice of the 2-norm is often chosen for ease
of computation and to guarantee the closed form solution, however,
statistically we should consider other norms. The choice of $\phi$
and $\psi$ depends on  the desired properties of the solution,
i.e., respectively, the choices of $\phi$ and $\psi$ represent
the class the solution $(w,\gamma )^\intercal$ should belong to and the
noise statistic, which may not be known. Recall that we wish to
address three aspects of the classification algorithms typically
used.  Each aspect will be considered separately, however the final algorithm
incorporates all elements.

\subsection{Improving sparsity via $\ell^p$ minimization}
\label{weightlp} 
The first
aspect considered is the sparsity of the weights $w$ corresponding to the neurons. The weights
$w$ largely consist of insignificant coefficients, which one may desire
to ``weed'' out, or effectively remove.  We introduce sparsity by
developing the $p$-norm method which weeds out unnecessary  weights
and selects the responsible data for each class. For the $p$-norm method we choose
\begin{equation} \label{lp}
\psi(w,\gamma)=|w|^p_p+\frac{1}{2}|\gamma|_2^2 \quad\mbox{where}\quad
|w|_p^p=\sum_{i=1}^m\;|w_i|^p
\end{equation}
for $0<p\le 1$. 
The
$p$-norm minimization has the effect of selecting the
desired data. In fact, it can be shown that
$$
|w|_p \to \mbox{\# of nonzero elements of $w$ as $p \to 0^+$}.
$$
Thus, $\psi$ enhances sparsity in the solution $w$ as $p\to 0^+$. In
other words, the choice of the $p$-norm with $p\leq 1$  removes
weights corresponding to neurons which are not active with respect to the
established baseline.

One can develop the necessary optimality condition for \eqref{Gen} with \eqref{lp}
despite the fact that $|w_i|_p$ is not differentiable at $w_i=0$.
Specifically, let us consider the case when $p=1$.
For $p=1$ we have
$$\d |w|=\dfrac{w}{|w|}$$ at $w\neq 0$, however for $w=0$ we have the subdifferential 
$$\d |w|=[-1,1]$$
where the subdifferential of a functional $f:\X\to (-\infty
,\infty]$ at $x\in\X$ is defined as the set \begin{equation} \{x^*\in\X^*| f(z)\geq
f(x)+\< x^*,z-x\>, \forall z\in\X\}, \label{subdiff} \end{equation} which can be found in \cite{itokun,convex} for example.
 However,
for $p<1$ we have $\d| w|_p^p = \emptyset$ when $w=0$.  To remedy
this, for $\varepsilon \ll 1$, we take the approximation \begin{equation} \d_\varepsilon
| w|_p^p = \dfrac{p\;w}{\max(\varepsilon^{2-p},|w|^{2-p})}\approx \d|w|_p^p	
\label{} \end{equation} which approximates the formal derivative \begin{equation} \d | w|_p^p = \dfrac{w}{|w|^{2-p}}
\label{} \end{equation}
for any value of $w$.
One can see a depiction of the approximate derivative $\d_{\varepsilon}|w|$ in Figure \ref{subdiff} 
with $\varepsilon=10e^{-3}$ and $\varepsilon=10e^{-2}$.
\begin{figure}[!hh]
\begin{center}
\includegraphics[width=3.0in]{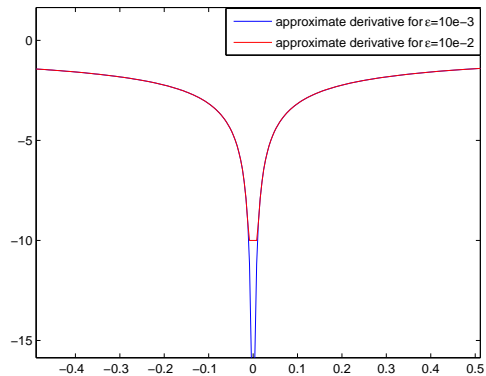}
\end{center}
\caption{Comparison of approximate subdifferentials for $\varepsilon=10e^{-3}$ and $\varepsilon=10e^{-2}$.}
\label{subdiff}
\end{figure}

Based on this, consider the iterative algorithm of the form 
\begin{equation}\label{iterate}
H^*Hw^{k+1}+\dfrac{\beta p}{\max(\varepsilon^{2-p},|w^k|^{2-p})}w^{k+1}=H^*e
\end{equation}
to find the minimizer of \eqref{Gen} with \eqref{lp}.
It can be shown that this sequence
converges to the minimizer of the appropriate cost functional, which is discussed in Appendix \ref{appenda}.

\subsection{Improving separability via inequality constraint}
\label{inequalityalg} 
We now present an algorithm which indirectly
utilizes the inequality $D(Aw-e\gamma)+y\ge e$ for the
optimization of $(w,\gamma)^\intercal$.  For this approach, we consider the
constrained minimization
\begin{eqnarray}
& \min\limits_{(w,\gamma)} & \phi(y)+\beta\psi(w,\gamma) \\
& \mbox{subject to } &D(Aw-e\gamma)+y= e \label{inequality}
\end{eqnarray}
where we design $\phi$ to incorporate the weighting for regions
where $(Aw-e\gamma)+y<e$ and $(Aw-e\gamma)+y>e$.  In this case,
the functional $\phi$ is taken to  be 
\begin{equation} \phi(y) = \dfrac{1}{2}\min(0,Hu-e)^2. \label{minphi} \end{equation} 
Here, we are attempting to
create more of a division between the two classes of data, so that
the data is classified more distinctly.
The choice \eqref{minphi} works well since if $y\ge 0$ we will have
$$(x\cdot w - \gamma)=1-y_i$$
which means that  $y_i$ is the degree of misclassification,
however, if $y_i\le 0$ we have
$$(x_i\cdot w-\gamma)=1-y_i\ge 1$$
which is desirable since the data is pushed farther away from 1.

\subsection{Weights for reducing bias}
Due to the large amount of data for periods of rest, the standard classification algorithm could have a bias towards identifying the
 ``rest'' state, hence reducing
the identification of an imagined movement.  The rest state corresponds to data for which $d_i=-1$.  Hence, it is reasonable
to consider incorporating different weights for the two cases $d_i=-1$ and $d_i=1$.  To reduce the bias towards identifying coefficients 
for which $d_i=-1$, we choose a parameter $\alpha\le 1$ for the weight corresponding to this data.  
In that way, for $d_i=-1,$ we weight the data by selecting 
a parameter $\alpha$ and we define the norm
$$
\phi(y) = \dfrac{1}{2}\| S^{1/2}(Hu-e)\|^2
$$ 
where the matrix $S$ is defined by
$$\left\{\begin{array}{l l}
S_{ii}=1\quad & \mbox{if } d_i=1 \\
S_{ii}=\alpha\quad & \mbox{if } d_i=-1
\end{array}\right.$$
in order to account for the bias.

\subsection{Algorithm}
We now provide the details of the numerical implementation of the methods discussed in Section \ref{newapproach} which contains
all the desired properties discussed.
For all algorithms presented in this section, the Tikhonov regularization parameter is $\beta$ and we define the matrix
\begin{equation}
H=D[A \,\,-e].
\label{}
\end{equation}

We develop an iterative algorithm which incorporates all aspects given in this section, based on
the iterative method \eqref{iterate}.
In summary, the iterative method for computing the
classifier $(w,\gamma)^\intercal$ is given by

\begin{equation}
\begin{array}{l}
\left(H^\intercal S\Gamma^k H+\beta\left(\begin{array}{cc} T^k & 0\\0
&1\end{array}\right)\right)\left(\begin{array}{c} w^{k+1}\\
\gamma^{k+1}\end{array}\right)=H^\intercal S\Gamma^k e
\\ \\
\Gamma^k_{ii}=\max(0,1-d_k(x_k-\gamma)),\quad T^k_{jj}=
\dfrac{p}{\max(\epsilon^{2-p},|w^k_j|^{2-p})}
\end{array}
\label{ouralg}
\end{equation}
for some small $\epsilon >0$.
Each step involves forming the diagonal matrices $\Gamma$, $T$
and solving the linear equation for $(w^{k+1},\gamma^{k+1})^\intercal$. 
For the application considered in this paper, convergence is achieved with a relatively low number
of iterations(3 or 4 is reasonable), so that the complexity of the proposed algorithm nearly
equals the complexity of the PSVM.  Thus, the advantages of this approach are realized with only a small increase
in computational cost for our application.  
Note that \eqref{ouralg} is equivalent to the PSVM formulation if we set $S,T,\Gamma=I\in\mathbb{R}^{n\times m}$.

\section{Choice rules for the regularization parameter}
\label{choicerules}
One of the questions that arises is how to choose the parameter $\beta=\frac{1}{\nu}$ optimally, 
so that the best solution is obtained.  Naturally, the choice of $\beta$ will depend on 
the choices of the functionals $\phi$ and $\psi$.  We present here not only the well known choice rule due to Morozov, but also a choice rule developed in 
\cite{regparam,newchoice} which has fewer assumptions than the well known Morozov's discrepancy principle.

If one knows the noise level in the data, then a well known and useful choice rule is the Morozov's Discrepancy Principle.  This choice rule works quite well in the case 
of known noise levels.
For the optimal solution $u_\beta$, the Morozov discrepancy principle seeks $\beta>0$ such that
$$\phi(u_{\beta},y)\simeq\delta$$
where $\delta$ is the  noise level defined by
$$ |y-y^\delta|_\mathcal{Y}\le\delta$$
for exact data, $y$, and noisy data, $y^\delta$.  Here, $\delta$ can be thought of as the performance level of the optimal solution.

It should be noted, that the approximate solutions $u_\beta^\delta\to u$ in the case that $\delta\to 0$.
Though the approximations converge, for problems such as classification, other choice rules may provide better performance. 
This has led to extensive research for the goal of developing choice rules which are able to automatically tune the parameter, $
\beta$ for the particular problem.

We now present a choice rule for automatically selecting the parameter $\beta$ based on the choices of $\phi,\psi$.
Consider maximizing the conditional density $p((u,\tau,\lambda)|y)\sim p(y|(u,\tau,\lambda))p(u,\tau,\lambda)$ where $(\tau,\lambda)$ are 
density functions for $\phi,\psi$, respectively, both having Gamma distribution.  The balancing principle is derived
from the Bayesian inference \cite{augtik}
$$
\min\limits_{(u,\tau,\lambda)} \tau\phi(u,y)+\lambda\psi(u)+\tilde{\beta}_0\lambda-\tilde{\alpha_0}\ln\lambda +\tilde{\beta}_1\tau-\tilde{\alpha}_1\ln\tau.
$$
Letting $\beta=\frac{\lambda}{\tau}$, the necessary optimality is given by
$$
\begin{array}{l}
u_\beta =\arg\,\min\limits_u\{\phi(u,y)+\beta\psi(u)\}\\
\\
\lambda =\dfrac{\tilde{\alpha}_0}{\psi(u_\beta)+\tilde{\beta}_0}\\
\\
\tau =\dfrac{\tilde{\alpha}_1}{\phi(u_\beta)+\tilde{\beta}_1}
\end{array}
$$
or
\begin{align}
u_\beta & = \arg\,\min\limits_u\{\phi(u,y)+\beta\psi(u)\}\\
& & \notag \\
\beta & = \dfrac{1}{\mu}\dfrac{\phi(u_\beta)+\tilde{\beta}_0}{\psi(u_\beta)+\tilde{\beta}_1},\quad \mu=\dfrac{\tilde{\alpha}_1}{\tilde{\alpha}_0.}\label{balanceprinc}
\end{align}
A variational formulation of the balance principle \eqref{balanceprinc} is developed in \cite{regparam,newchoice}.  
The theoretical justification for the balance principle follows from the variational formulation. Good choice rules
 for regularization parameters become increasingly important as the data becomes increasingly noisy and, likewise,
when the problem is highly ill-posed.  In the case that the statistic $\mu$ is unknown, an algorithm for estimating $\mu$
is discussed in \cite{regparam,newchoice}.

The natural choice for updating the parameter $\beta$ would be the fixed point iteration
\begin{enumerate}
\item Set $k=0$ and choose $\beta_0$
\item \label{step2} Solve for $u_k$ 
$$
\arg\,\min\limits_u \{\phi(u)+\beta_k\psi(u)\}
$$
\item Update the regularization parameter $\beta_{k+1}$ by
$$
\beta_{k+1}=\dfrac{1}{\mu}\dfrac{\phi(u_k)}{\psi(u_k)}
$$
\item Check the stopping criterion.  If convergence is not met, set $k=k+1$ and repeat from Step \ref{step2}.
\end{enumerate}

\section{Numerical Results}
\label{numres}
\begin{figure}[!hh]
\begin{center}
\includegraphics[width=2.75in]{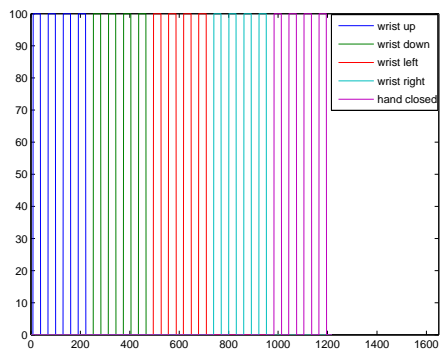}
\end{center}
\caption{Cues for the first patient.}
\label{cues2}
\end{figure}

In this section, we provide numerical results from the application of our approach to the neural classification problem outlined in the Introduction.
We first provide an outline of how the experiment is conducted, as it will be important for deciphering the results provided.
The experiment proceeds by asking a person think about a specific movement such as wrist flexion, elbow extension, or closing the hand.  The data obtained is
neural firing rate data, meaning that this is a time dynamical data classification problem.  It is desired to
determine the correlation between imagining a movement and the neural response to imagining the movement.  The experiment consists of periods of rest and periods
where the person is cued to think about a certain movement.  For this particular experiment, there are five wrist movements consisting of wrist extension,
wrist flexion, wrist radial deviation, wrist ulnar deviation, and closing hand.  The patient is cued to imagine one of the movements consecutively, 
with periods of rest between each cue.  The person is then cued for another of the five movements consecutively, again with periods of rest in between each cue.
For instance, an example set of cues for the movement of wrist up can be seen in Figure 
\ref{cues2}.  
The period of time after the cues for wrist up consists of both rest and cues 
for other movements.  This is the nature of how the experiment is conducted.  Data is collected for a specified interval of time, which consists of periods of rest and 
periods of cues for each movement.  The goal of the data classification is to sharply separate the data for each movement.  For example, we must separate the data 
corresponding to wrist ``up'' (extension) from the data corresponding to both the rest periods and the periods for other movements.  Thus, the results given here are only for one 
particular movement, however the method produces similar results for the other five movements.

We also develop and implement a new performance measure utilizing the force
\begin{equation}
F=Aw-\gamma.
\label{}
\end{equation}
Often, the performance measure is taken to be
\begin{equation}
\operatorname{sgn}(F)=\left\{\begin{array}{r l}
1\quad & \mbox{if } F>0\\
0\quad & \mbox{if } F=0\\
-1\quad & \mbox{if } F<0
\end{array}\right.
\label{}
\end{equation}
pointwise in time,
which could have a tendency to include outliers.  Instead, we consider other options for determining the performance.
One such option is to sum the force over a time interval. That is we define a performance measure $P$ by
$$P_k=\sum\limits_{i=k-n}^{k+n}F_i$$  for some $n$ (e.g. $n=5$).  One could also take the pointwise average as
$\frac{P}{2n+1}$, which is equivalent to $P$ as a performance measure since we are interested in the sign of that data.
We call $P$ the summed performance measure.
Another option would be to average the force over each time interval of cues and rest periods.  This performance 
measure is given by
\begin{equation}
\tilde{P}_k=\dfrac{1}{N}\sum\limits_{i=1}^{N}F_i
\label{}
\end{equation}
\begin{figure}[hbt]
\centering
\subfloat[][]{\label{perf-a}\includegraphics[width=2.25in]{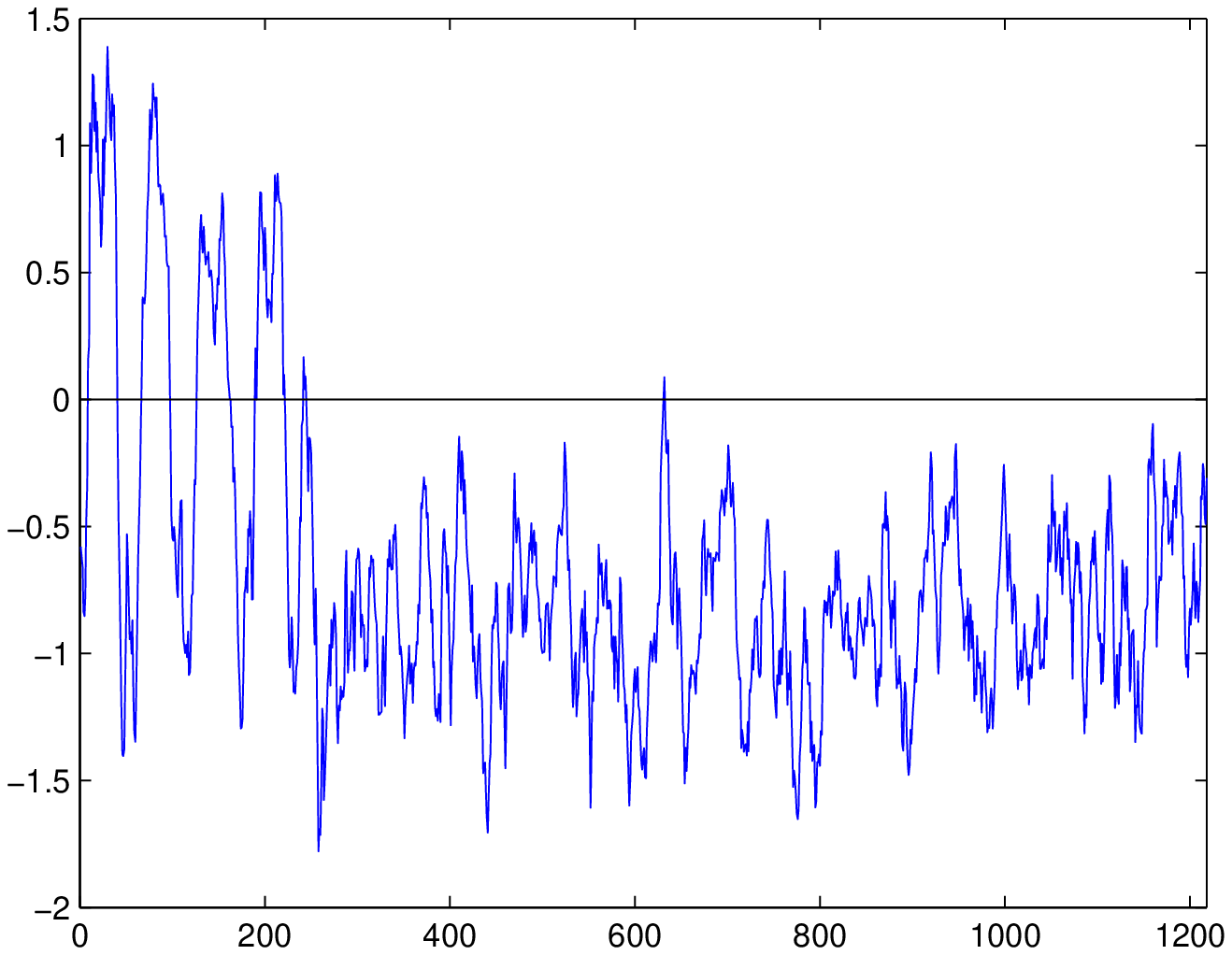}}\quad
\subfloat[][]{\label{perf-b}\includegraphics[width=2.25in]{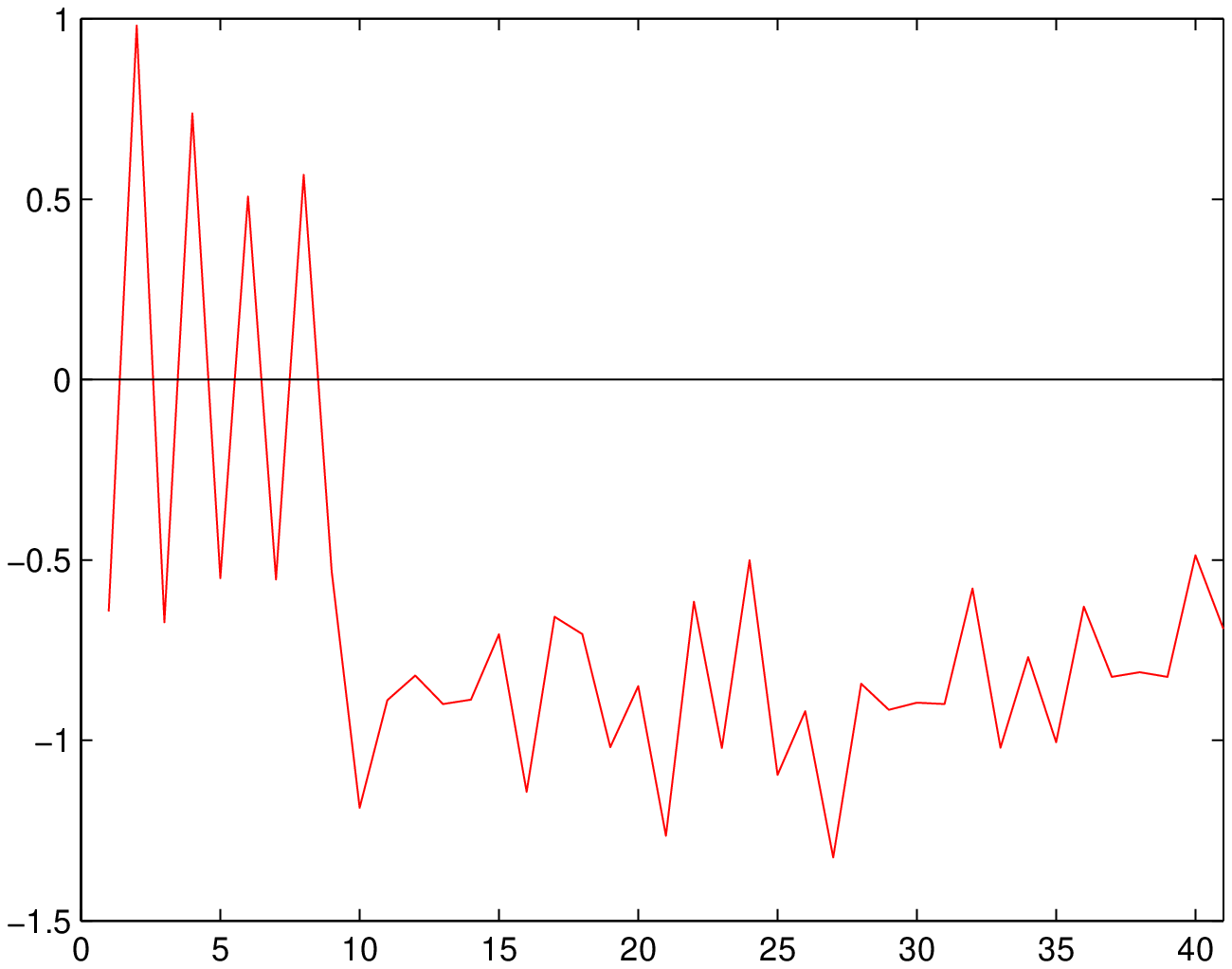} }\\
\subfloat[][]{\label{perf-c}\includegraphics[width=2.25in]{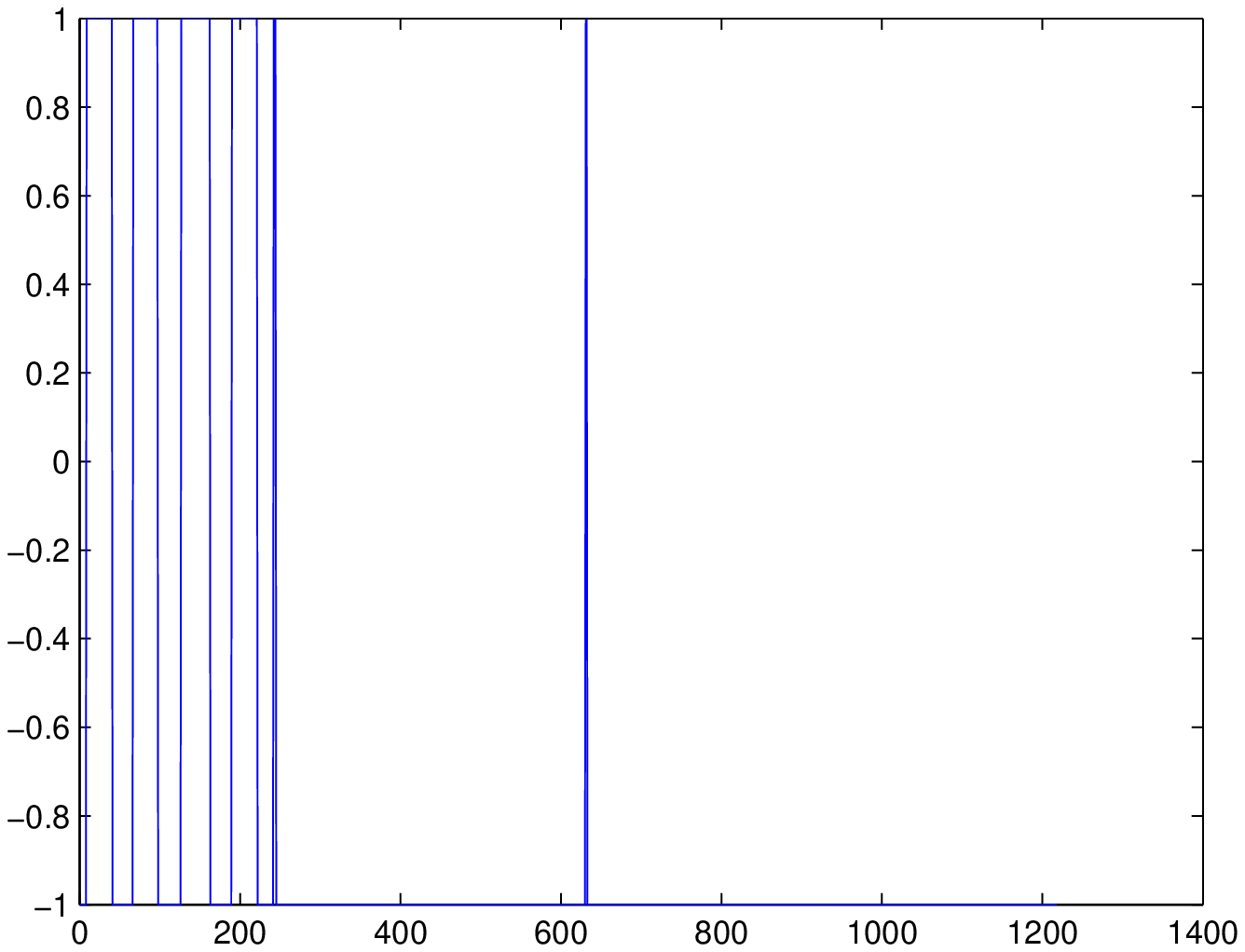}}\quad
\subfloat[][]{\label{perf-d}\includegraphics[width=2.25in]{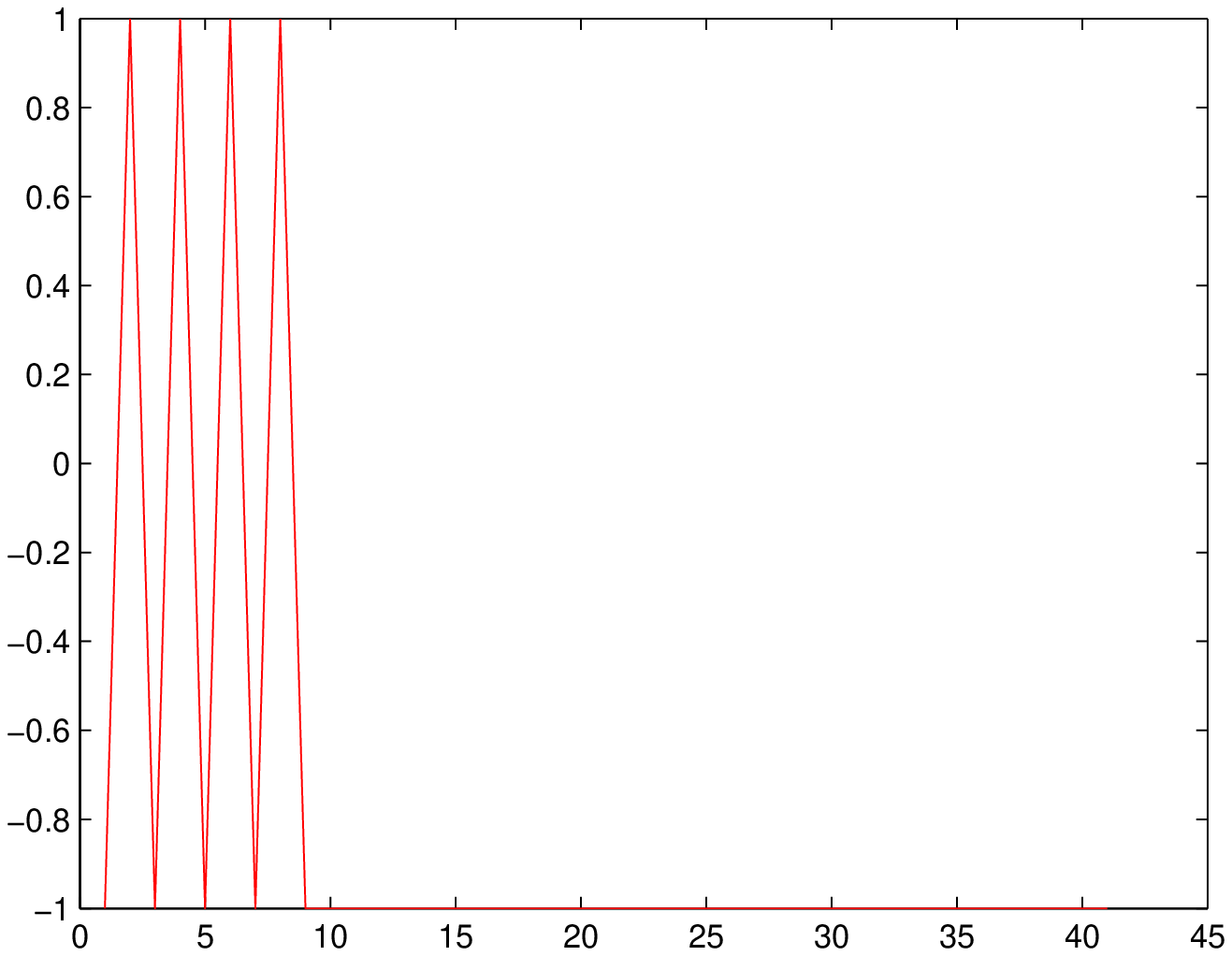} }
\caption{\subref{perf-a} An example of the force $F$; \subref{perf-b} An example of the averaged force;
 \subref{perf-c} The sign function of the original force; \subref{perf-d} The sign function of the averaged force.}
 \label{perf}
\end{figure}
where $(t_{1}\ldots t_{N})$ is the interval of time for the cue or rest period, and $N$ is the number of points.
Note that we are taking the average over each interval for which the sign of the surrounding data is the same.
Both $P$ and $\tilde{P}$ improve the detection since outliers are likely 
single points or small clusters.  The process of summing the forces can cause outliers to cancel out over the interval.
Likewise, the process of averaging over an interval decreases 
their effect.  For example, one can see in Figure 
\ref{perf} how $\operatorname{sign}(F)$ has an outlier while $\operatorname{sign}(\tilde{P})$ does not.
Analyzing the force yields the detection of the responsible neurons for a given movement, and hence the performance
of the algorithm.  

\begin{figure}
\centering
\subfloat[][]{\label{s0}\includegraphics[width=2.5in]{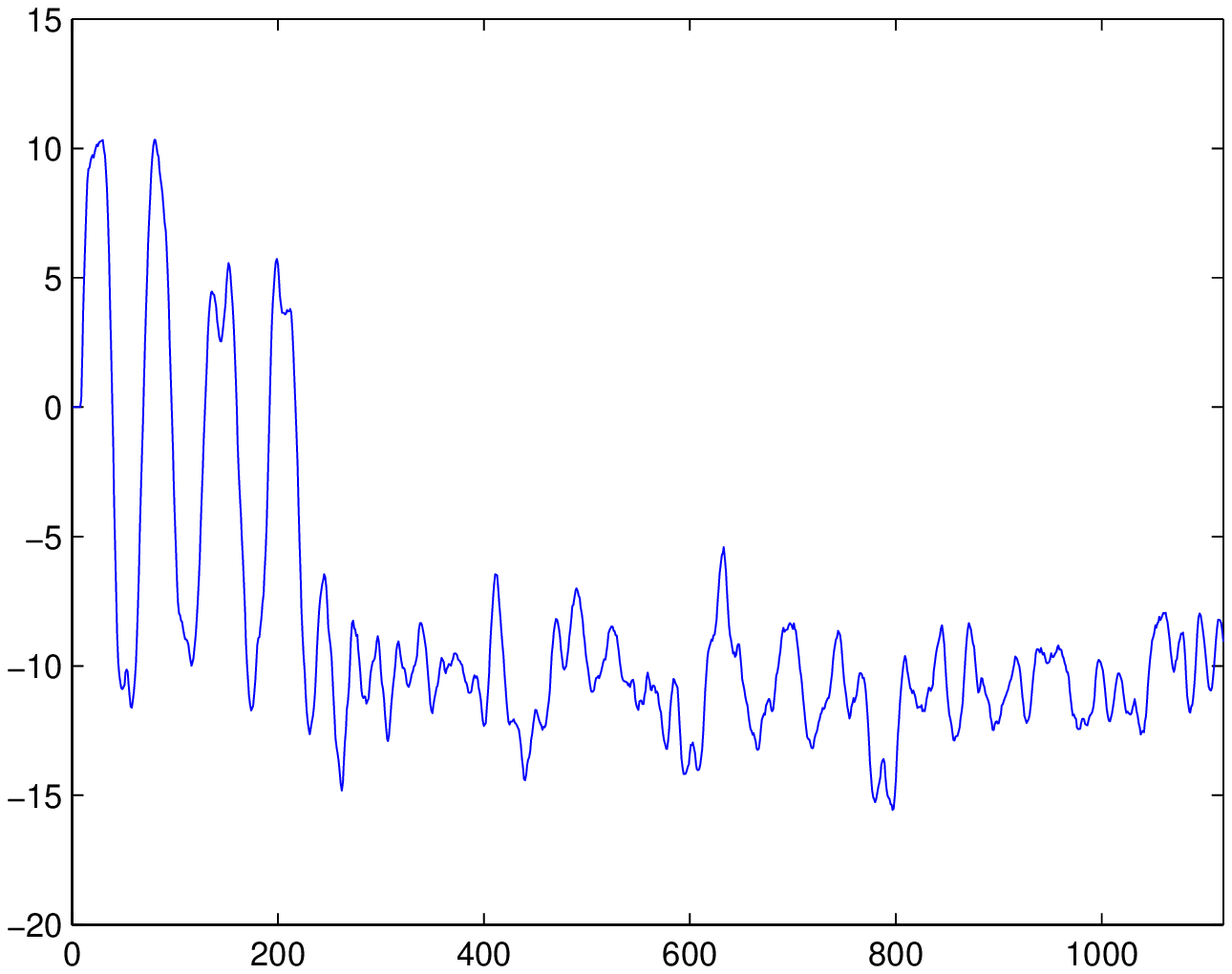}}\quad
\subfloat[][]{\label{M0}\includegraphics[width=2.5in]{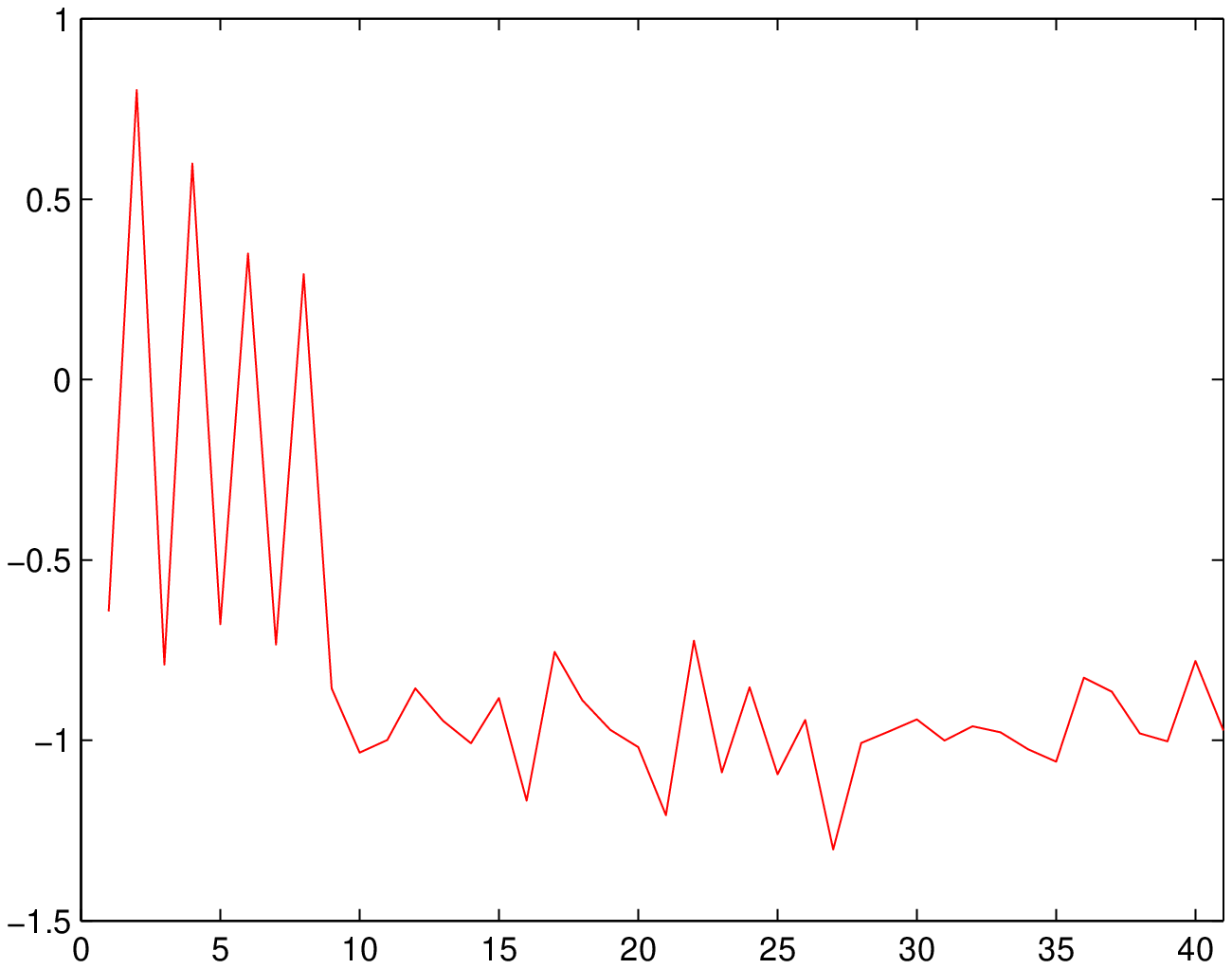} }
\caption[]{\subref{s0} Summed forced example \subref{M0} Averaged force example.}
 \label{s0M0}
\end{figure}

The two performance measures $P$ and $\tilde{P}$ are both likely good choices, however $\tilde{P}$ may distinguish more clearly between
the data classes, as can be seen in Figure 
\ref{s0M0}.

\subsection{Stroke Patient}
\begin{figure}[!hh]
\begin{center}
\includegraphics[width=2.75in]{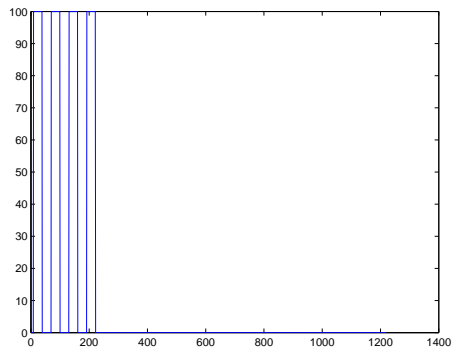}
\end{center}
\caption{Cues for the wrist up action.}
\label{cues}
\end{figure}

We now provide the numerical results for one of the patients who was left paralyzed after a stroke.  
As far as the two data sets used for this paper, the stroke patient's 
data is the easier of the two sets.  We applied the modifications to the original PSVM separately and some in 
conjunction with one another.  We provide comparisons of the force itself, the performance measure and we also include a heat map of the neural activity, which can be used 
for visualizing the active neurons.  The cues for the action of 'wrist up' can be seen in Figure 
\ref{cues} for this particular patient.

For clarity of the presentation, we separately consider the three improvements of our formulation, namely the increased sparsity, reduced bias towards classifying the at 
rest data and the robustness of separating the classes.  We provide results for each of the three separately, however the best implementation takes advantage of all three 
aspects. It should also be noted that all results for this patient are for the wrist up movement, unless otherwise stated. 

\begin{figure}[!hh]
\begin{center}
\includegraphics[width=2.75in]{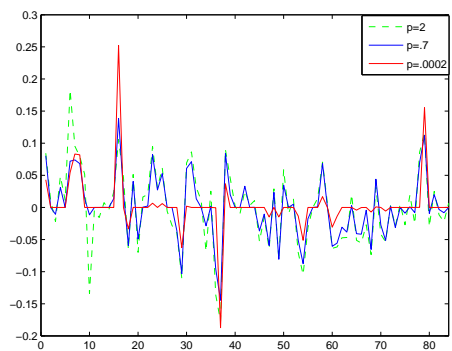}
\end{center}
\caption{Comparison of weights for different $\ell_p$ norms.}
\label{sparseweights}
\end{figure}

\begin{figure}[!hh]
\centering
\subfloat[][]{\label{l2vlp-a}\includegraphics[width=2.5in]{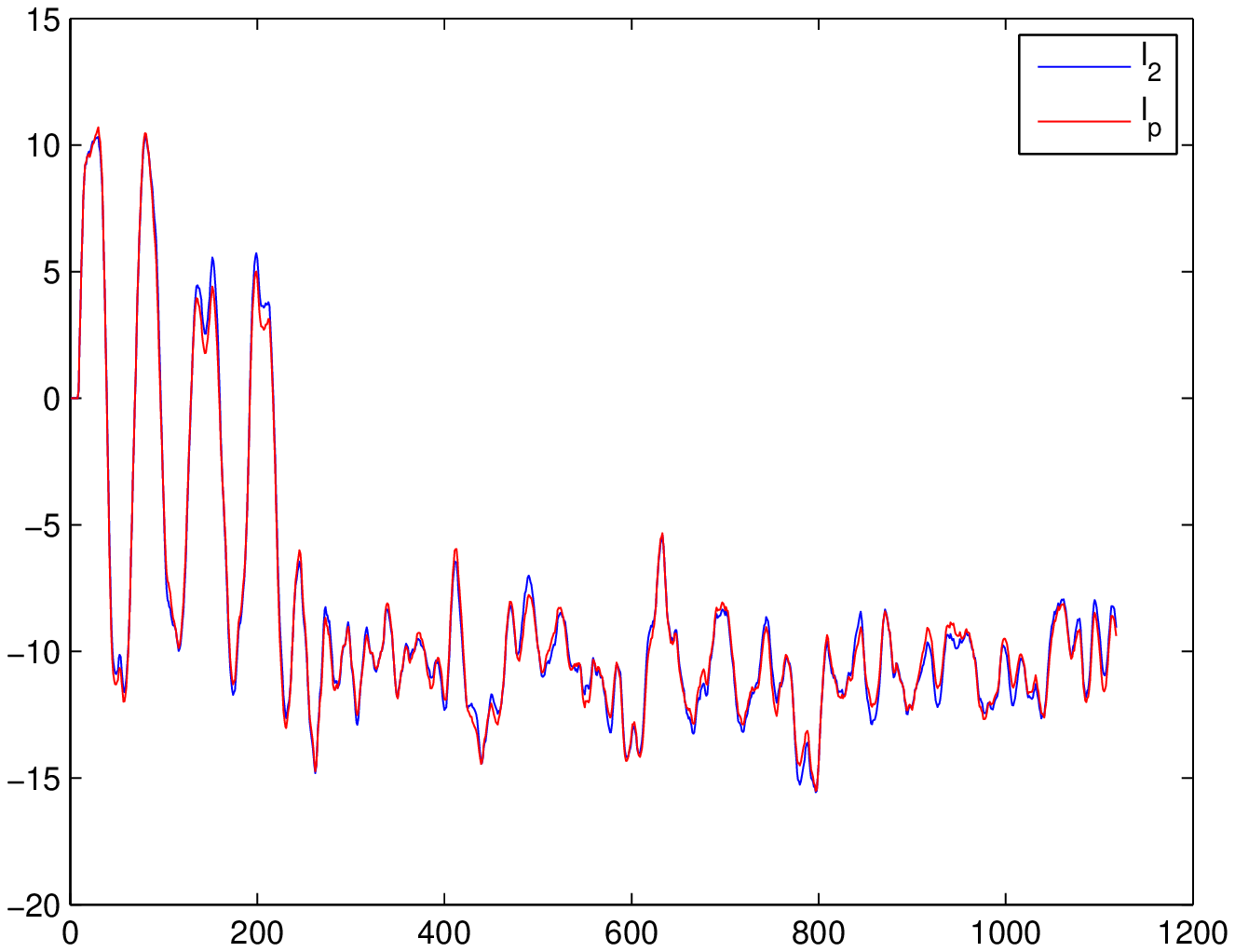}}\quad
\subfloat[][]{\label{l2vlp-b}\includegraphics[width=2.5in]{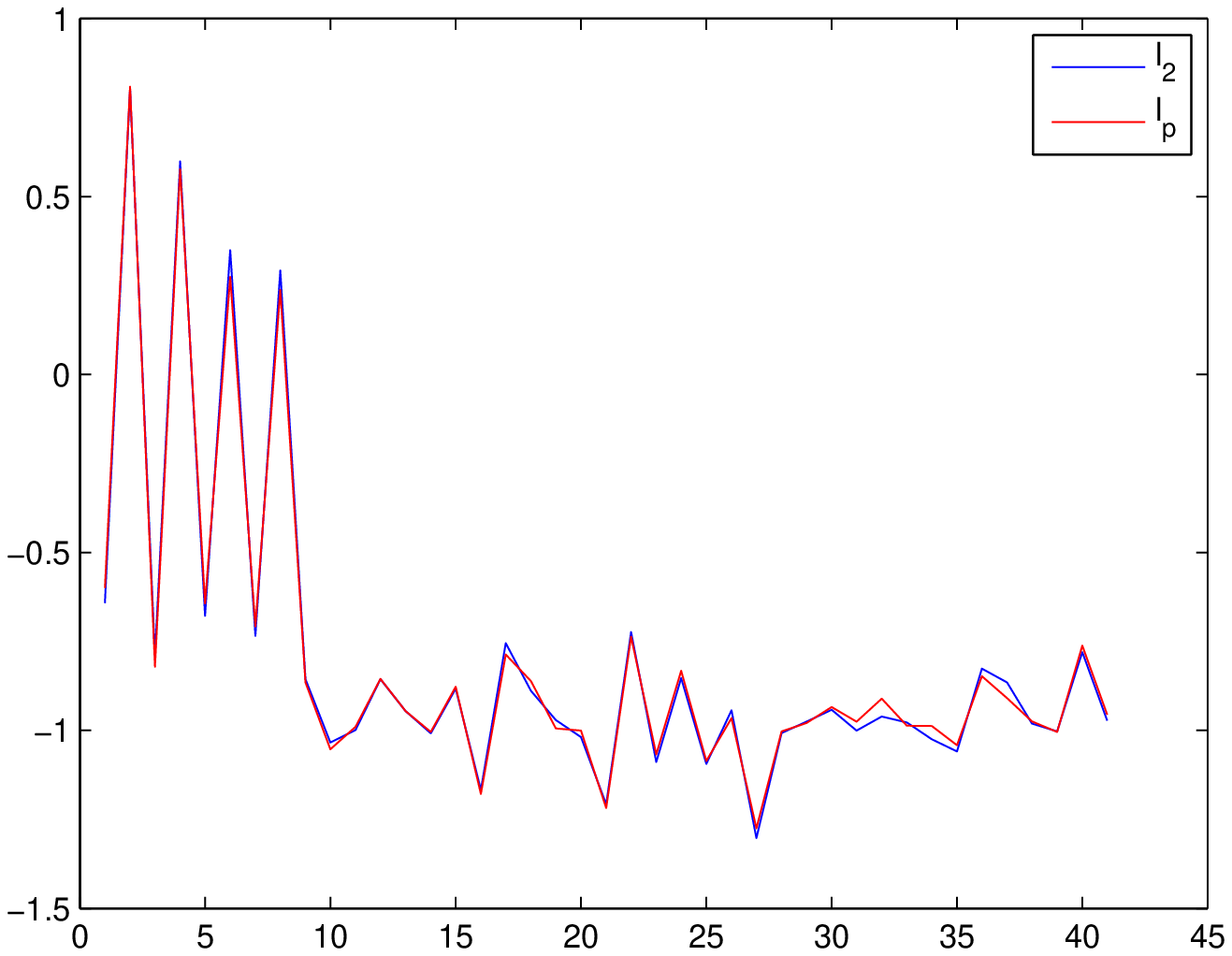} }
\caption[]{\subref{l2vlp-a} Comparison of $\ell_2$ versus $\ell_p, p=.2$ using summed force; 
\subref{l2vlp-b} Comparison of $\ell_2$ versus $\ell_p, p=.2$ using averaged force.}
 \label{l2vlp}
\end{figure}

By taking the $\ell_p$ minimization with $p=.2$ versus the $\ell_2$ minimization we are able to increase the number of coefficients such that
$|w_i|\le 1e-4$ from 1 to 18.  The increase in sparsity is depicted in Figure \ref{sparseweights} where one can see how taking successively 
smaller values of $p$ reduces the number of nonzero weights, $w$.  
We now illustrate how 
changing $\beta$ can improve the force.  By taking $\beta=.2$  one can see in Figure 
\ref{l2vlp} that the change in force by taking the $\ell_p$
 minimization versus the $\ell_2$ minimization is negligible, however an increase in sparsity is still realized.

\begin{figure}[!hh]
\centering
\subfloat[][]{\label{alph-a}\includegraphics[width=2.5in]{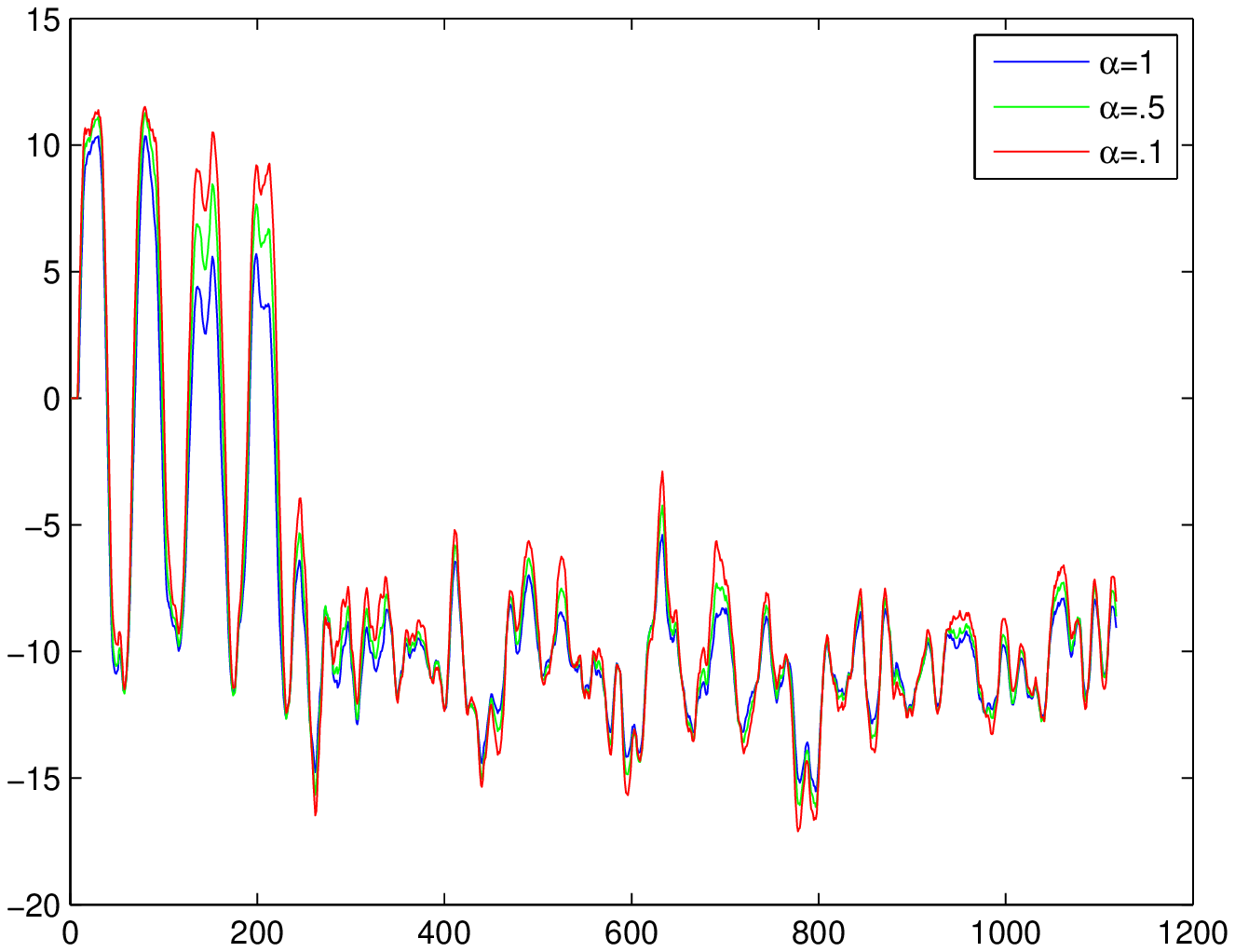}}\quad
\subfloat[][]{\label{alph-b}\includegraphics[width=2.5in]{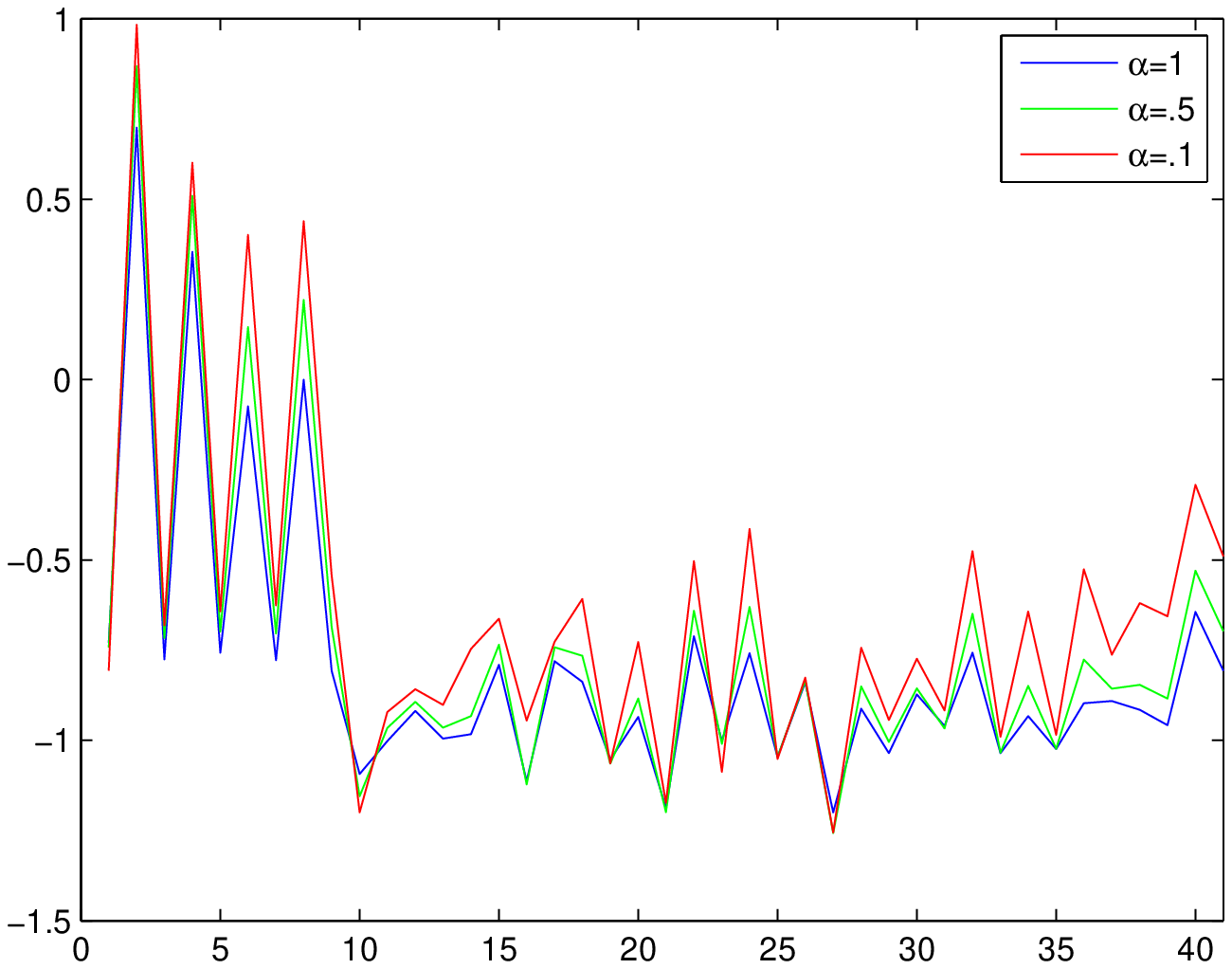} }
\caption[]{\subref{alph-a} Summed force for different values of $\alpha$; 
\subref{alph-b} Averaged force for different values of $\alpha$.}
 \label{alph}
\end{figure}

Now, let us consider the incorporation of a parameter $\alpha$ in order to reduce the bias towards classifying at rest data(or data for other movements).
If the value of $\alpha$ is chosen appropriately, the magnitude of the force increases where $F>0$, while the increase in magnitude is minimal for 
$F<0$.  This is illustrated in Figure \ref{alph}.
  
We now present results for improving the robustness of separating the classes, via the choice \eqref{minphi}.  In this case, an increase in magnitude
for $F>0$ and a decrease in magnitude for $F<0$ is realized.  This is ideal, since it is desired to completely separate the classes.  Note that the forces
have been normalized to present a fair comparison.

\begin{figure}[!hh]
\centering
\subfloat[][]{\label{yym-a}\includegraphics[width=2.5in]{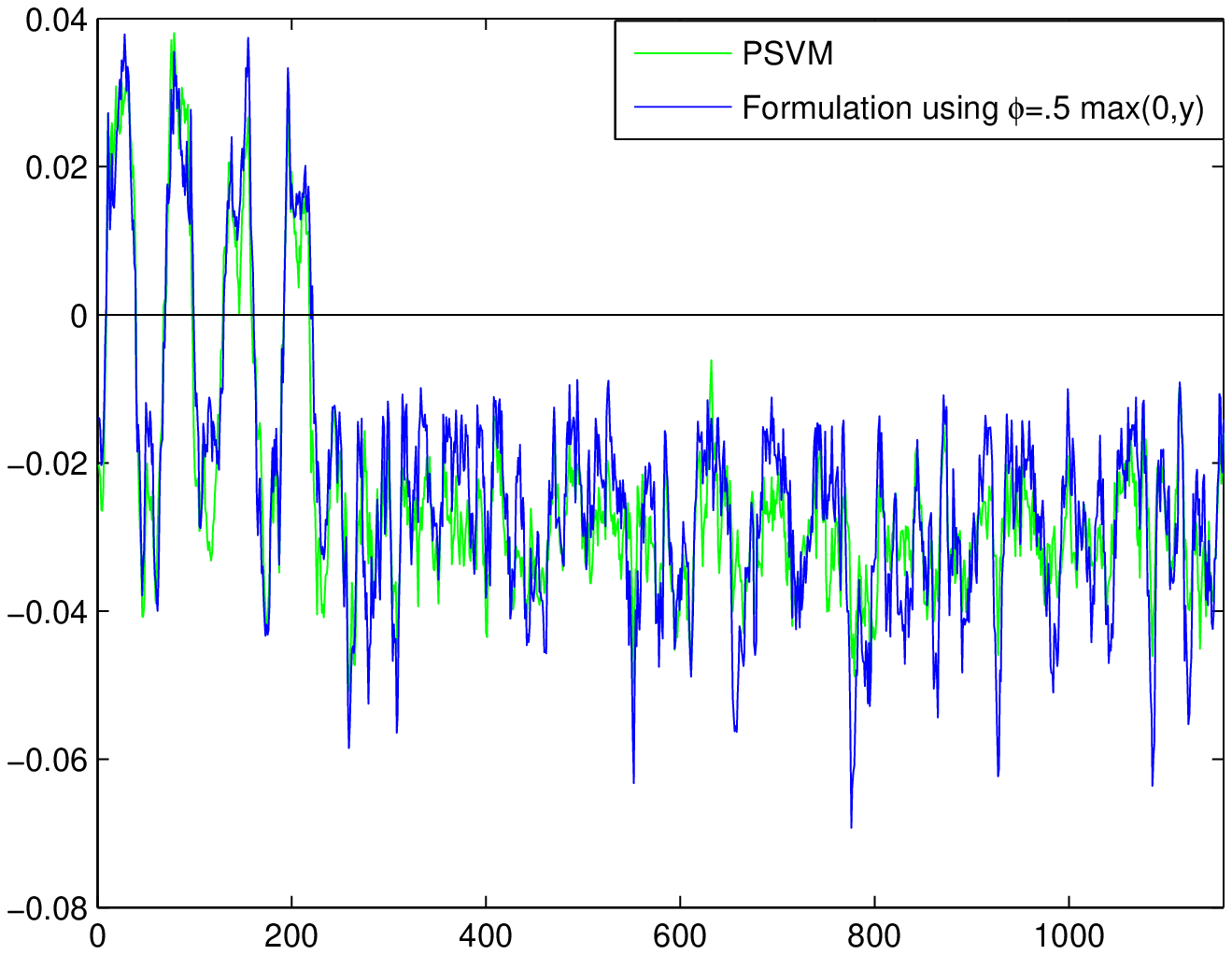}}\quad
\subfloat[][]{\label{yym-b}\includegraphics[width=2.5in]{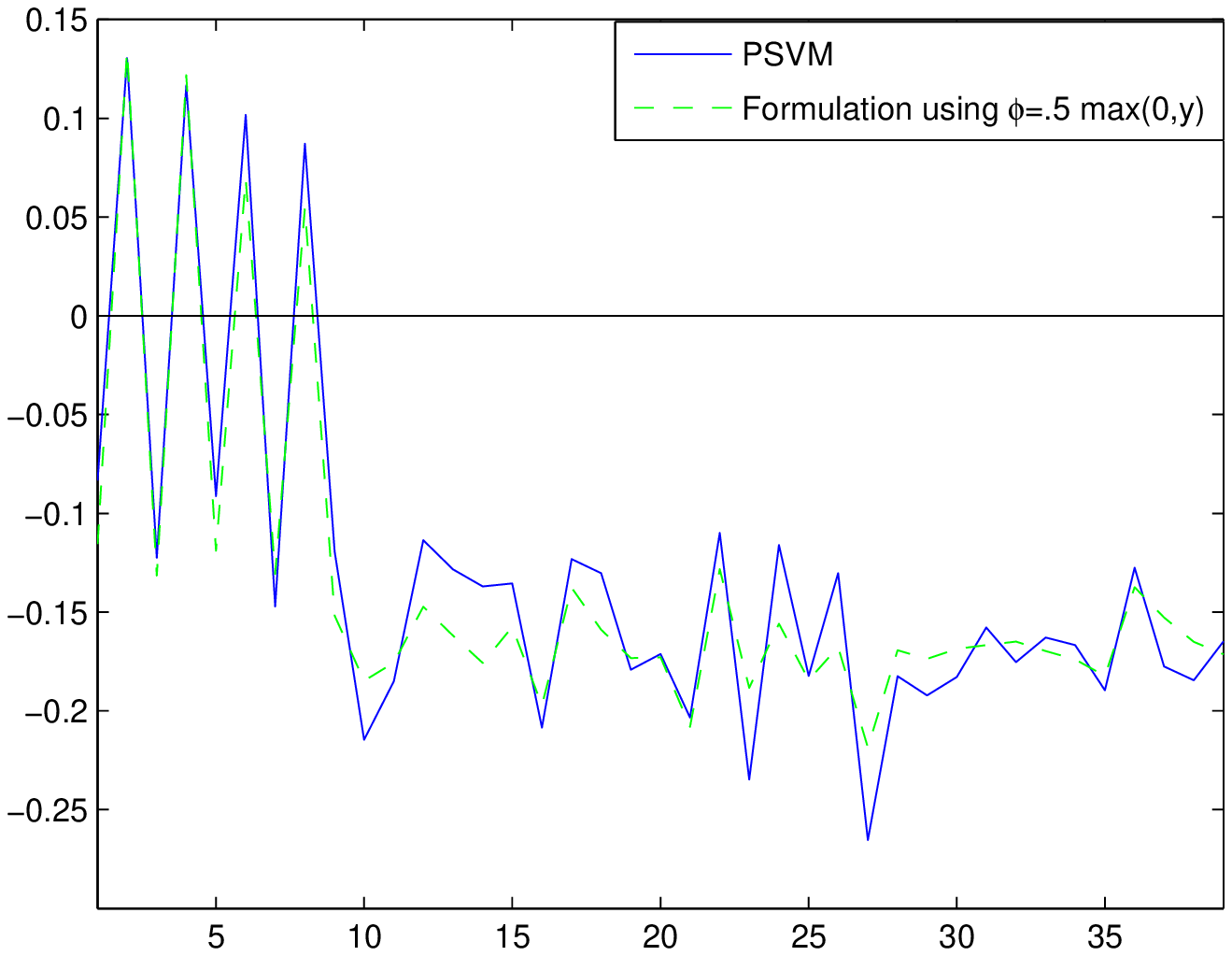} }
\caption[]{\subref{yym-a} Summed force after minimization of $\phi(y^+)$ (normalized) for $\alpha =1, p=1$; 
\subref{yym-b} Averaged force after minimization of $\phi(y^+)$ (normalized) for $\alpha =1, p=1$.}
 \label{yym}
\end{figure}

\begin{figure}[!tt]
\begin{center}
\includegraphics[width=3.5in]{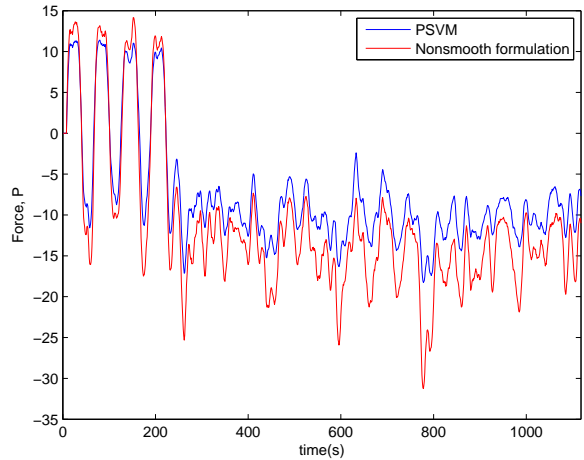}
\end{center}
\caption{Comparison of PSVM and nonsmooth formulation for $\alpha =.5, p=.2$.}
\label{psvmvsnon}
\end{figure}

Next, we display a comparison of the summed force, $P$, for the PSVM formulation versus the nonsmooth formulation developed in this paper, 
in Figure 
\ref{psvmvsnon}.  Note that the PSVM is the same as \eqref{Our} with $S,T,\Gamma=I\in\mathbb{R}^{n\times m}$.  As can be seen, the separation of the considered movement 
from all other movements is improved.
  
\begin{figure}
\centering
\subfloat[][]{\label{wristtest-a}\includegraphics[width=2.5in]{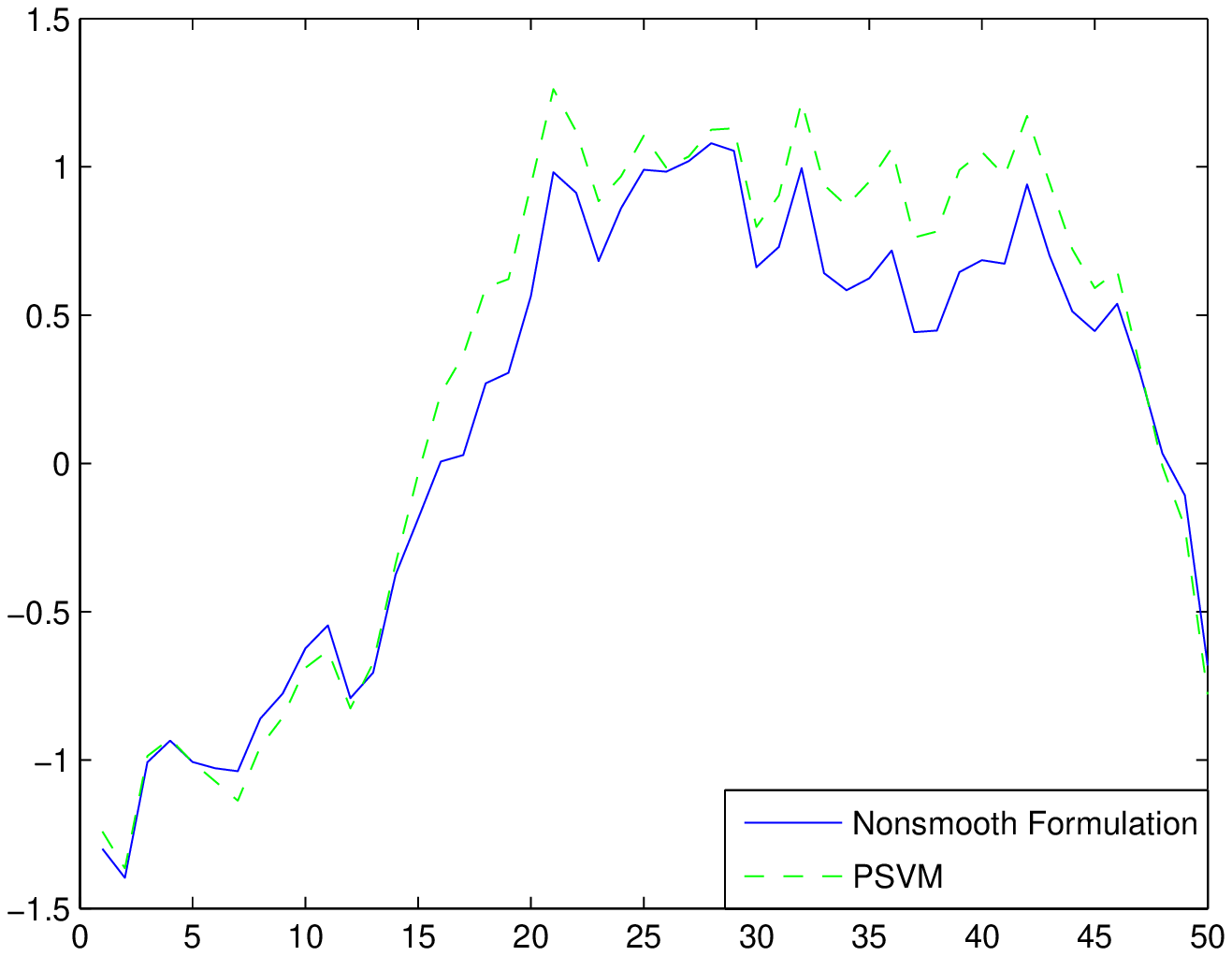}}\quad
\subfloat[][]{\label{wristtest-b}\includegraphics[width=2.5in]{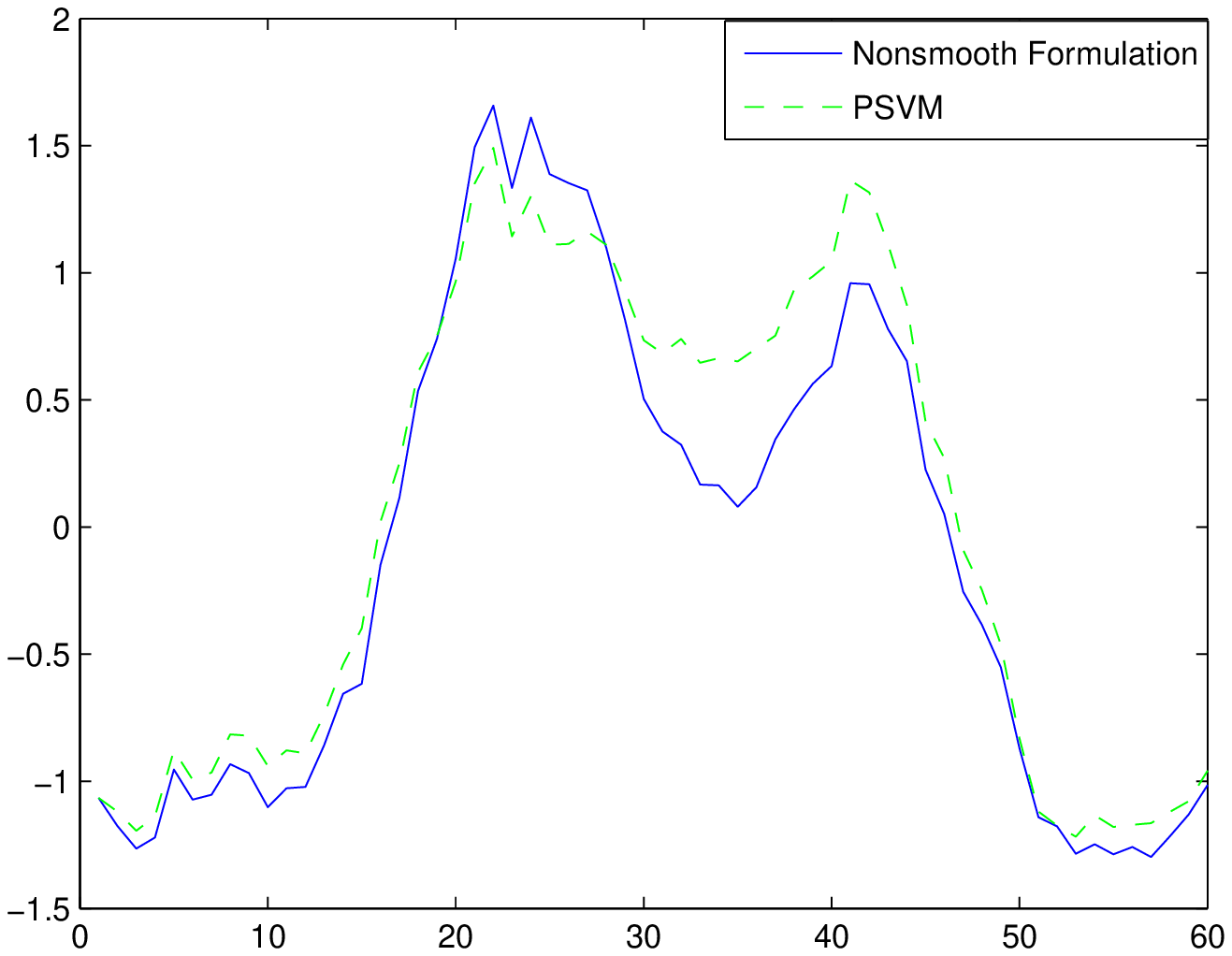} }
\caption[]{\subref{wristtest-a} Test of nonsmooth formulation for wrist down for $\alpha =.1, p=2$;
 \subref{wristtest-b} Test of nonsmooth formulation for wrist right for $\alpha =.1, p=2$.}
 \label{wristtest}
\end{figure}

As a test of the classifier, we train on the first three cues for a specified movement and test the classifier
on the fourth cue.  The test results can be seen in Figure \ref{wristtest} for two movements.  
As one can see, the classification is nearly equal for both the PSVM and our nonsmooth formulation.

\vspace{2cm}
\subsection{Lou Gehrig's Disease Patient}
We now provide results for a patient affected by Amyotrophic Lateral Sclerosis(ALS), or Lou Gehrig's Disease.  The majority of the results are similar to those obtained 
for the stroke patient, with the exception of one outlier in the results.  Several explanations should be considered as to why the outlier appears in the results.  One 
potential explanation is that the patient suffered from performance degradation due to the increased length of this particular experiment.  Secondly, the nature of this 
patients physical disability is different from that of the first patient, so there may be difficulty with respect to the natur of the disease itself.
Even with this outlier, the results are improved, 
considering that the outlier is a single spike in the results and can likely be ignored when analyzing the data.  This particular data set consists of two periods of cues, as 
opposed to the one period for the first data set.  Being the more difficult of the two data sets, this truly tests the power of our nonsmooth formulation.

\begin{figure}[!hh]
\begin{center}
\includegraphics[width=2.75in]{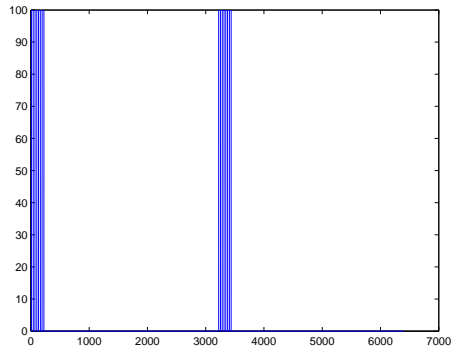}
\end{center}
\caption{Cues for the wrist up action.}
\label{cuesals}
\end{figure}

\begin{figure}[!hh]
\centering
\subfloat[][]{\label{l2vlpals-a}\includegraphics[width=2.5in]{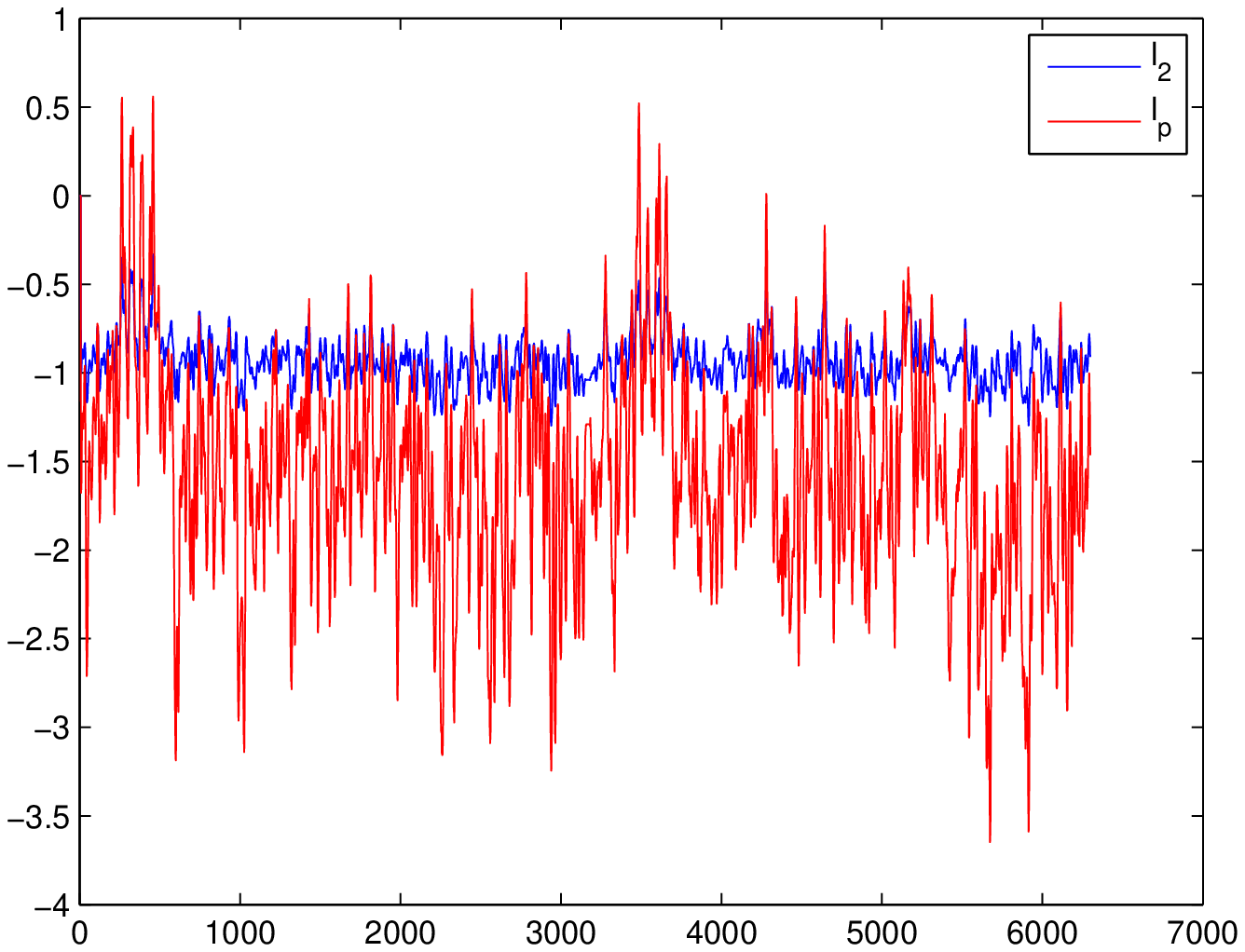}}\quad
\subfloat[][]{\label{l2vlpals-b}\includegraphics[width=2.5in]{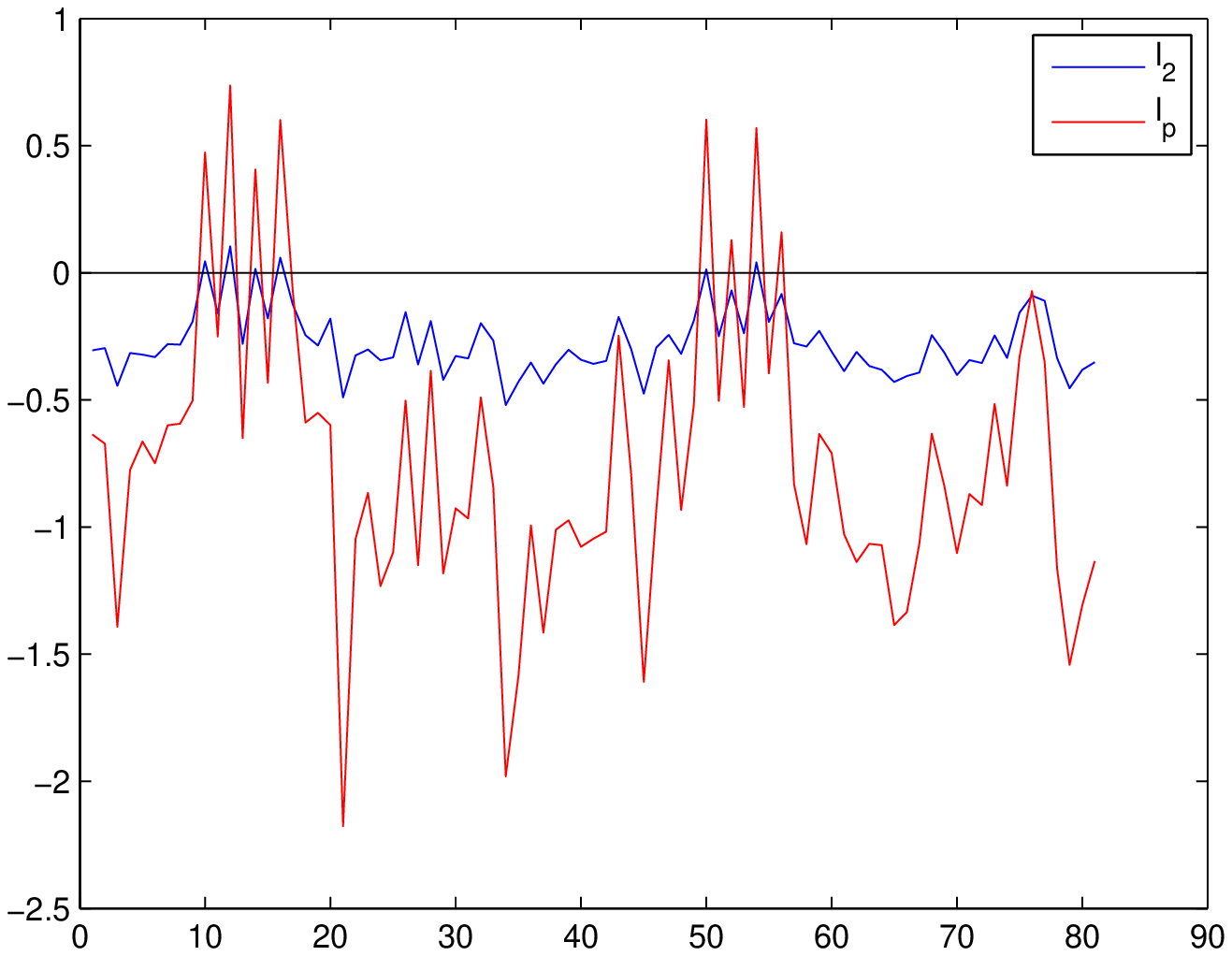} }
\caption[]{\subref{l2vlpals-a} Comparison of $\ell_2$ versus $\ell_p$ using summed force for $p=.2$;
 \subref{l2vlpals-b} Comparison of $\ell_2$ versus $\ell_p$ using averaged force for $p=.2$.}
 \label{l2vlpals}
\end{figure}

\begin{figure}[!hh]
\centering
\subfloat[][]{\label{alphals-a}\includegraphics[width=2.5in]{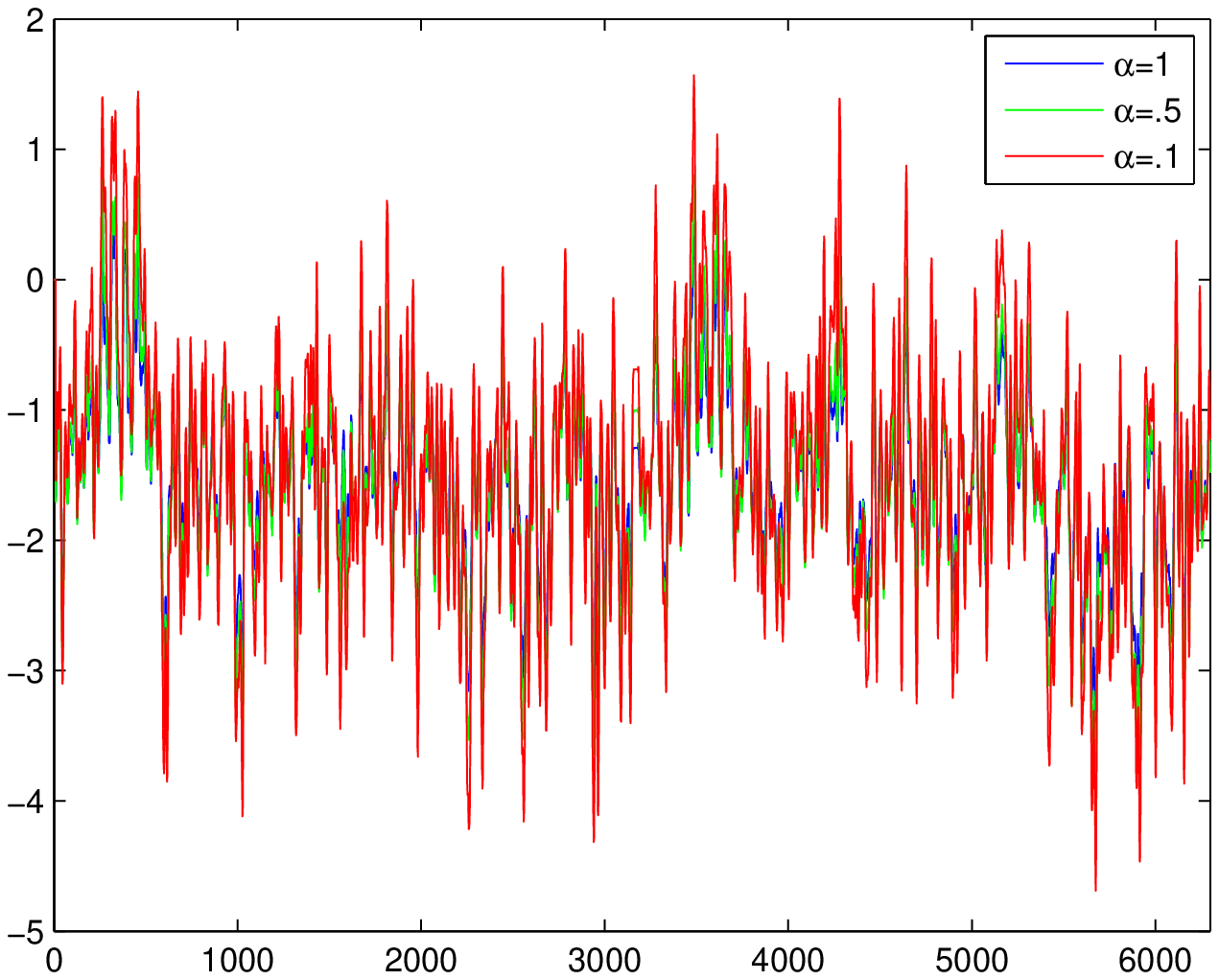}}\quad
\subfloat[][]{\label{alphals-b}\includegraphics[width=2.5in]{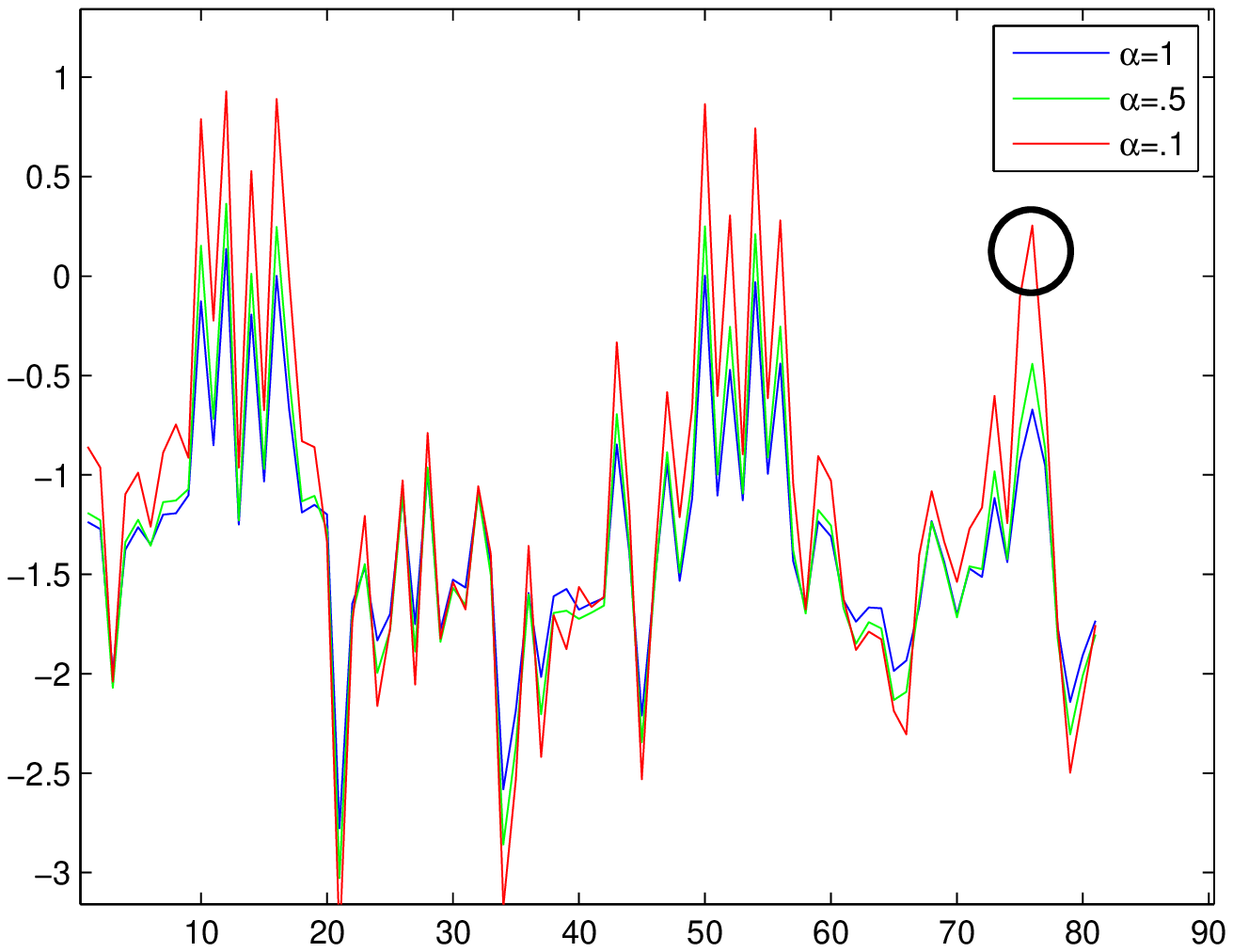} }
\caption[]{\subref{alphals-a} Summed force for different values of $\alpha$  \subref{alphals-b} Averaged force for different values of $\alpha$ .}
 \label{alphals}
\end{figure}

One can see similar results for this patient in Figures 
\ref{l2vlpals} and 
\ref{alphals}, where one can see the circled outlier in Figure 
\ref{alphals-b}.
This particular example shows how the tuning of the parameter $\alpha$ can affect the classification results.  
If $\alpha$ is chosen to be too small (e.g. $\alpha=.1$ for this data set)
then the appearance of outliers may increase.

\begin{figure}[!hh]
\centering
\subfloat[][]{\label{testals-a}\includegraphics[width=2.5in]{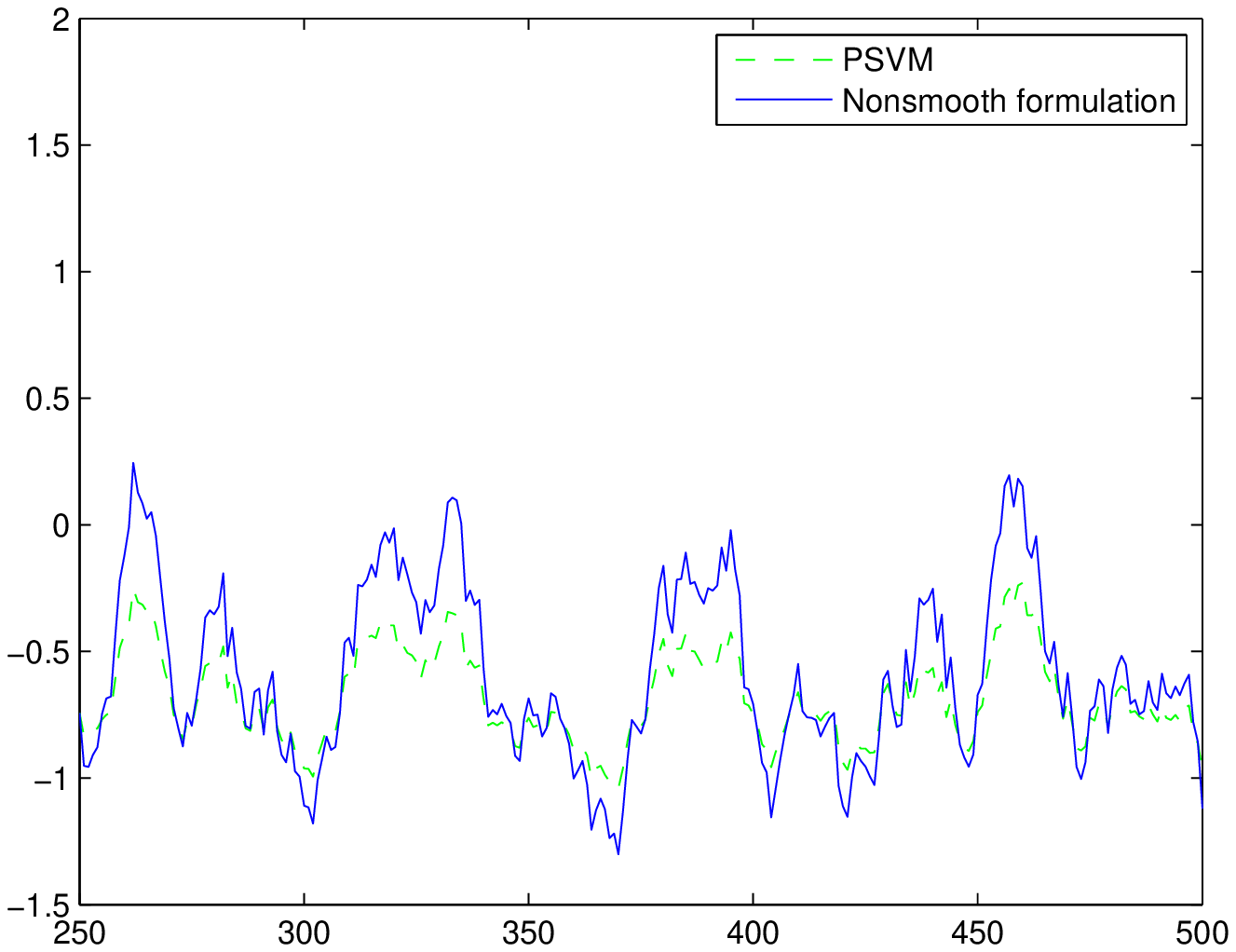}}\quad
\subfloat[][]{\label{testals-b}\includegraphics[width=2.5in]{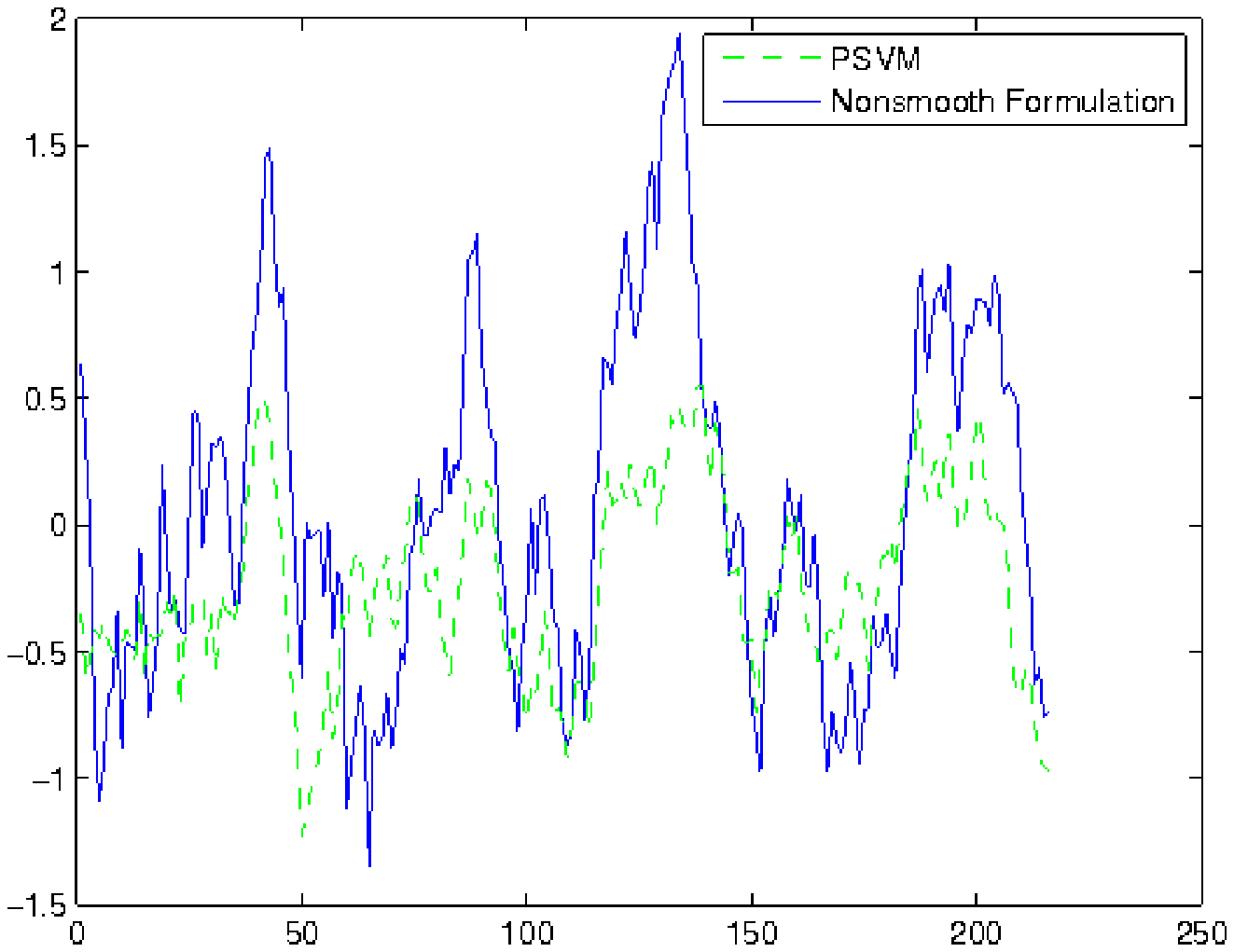} }
\caption[]{\subref{testals-a} Training corresponding to test of nonsmooth formulation for $\alpha =.1, p=2$;
\subref{testals-b} Test of nonsmooth formulation for wrist up for $\alpha =.1, p=2$.}
\label{testals}
\end{figure}

As a test of the classifier, we train on the first set of cues for a specified movement and test the classifier
on the second set of cues.  The test results can be seen in Figure \ref{testals} for one movement.  
As one can see, the classification improves by using the nonsmooth formulation versus
the PSVM.

\subsection{Identifying the Physical Location via Heat map}
\begin{figure}[!hh] 
\centering 
\subfloat[][]{\label{sparseheat-a}\includegraphics[width=2.5in]{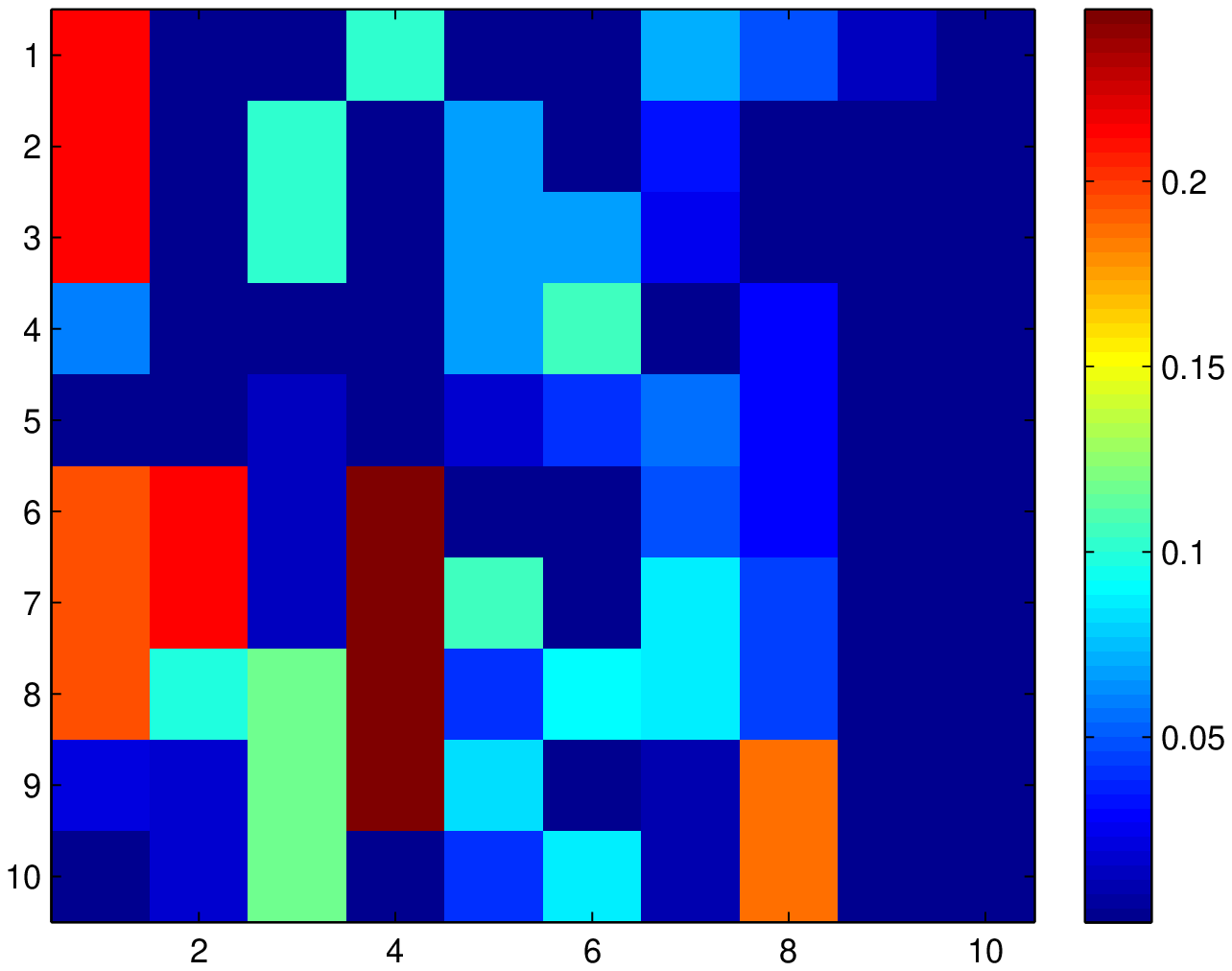}}\quad 
\subfloat[][]{\label{sparseheat-b}\includegraphics[width=2.5in]{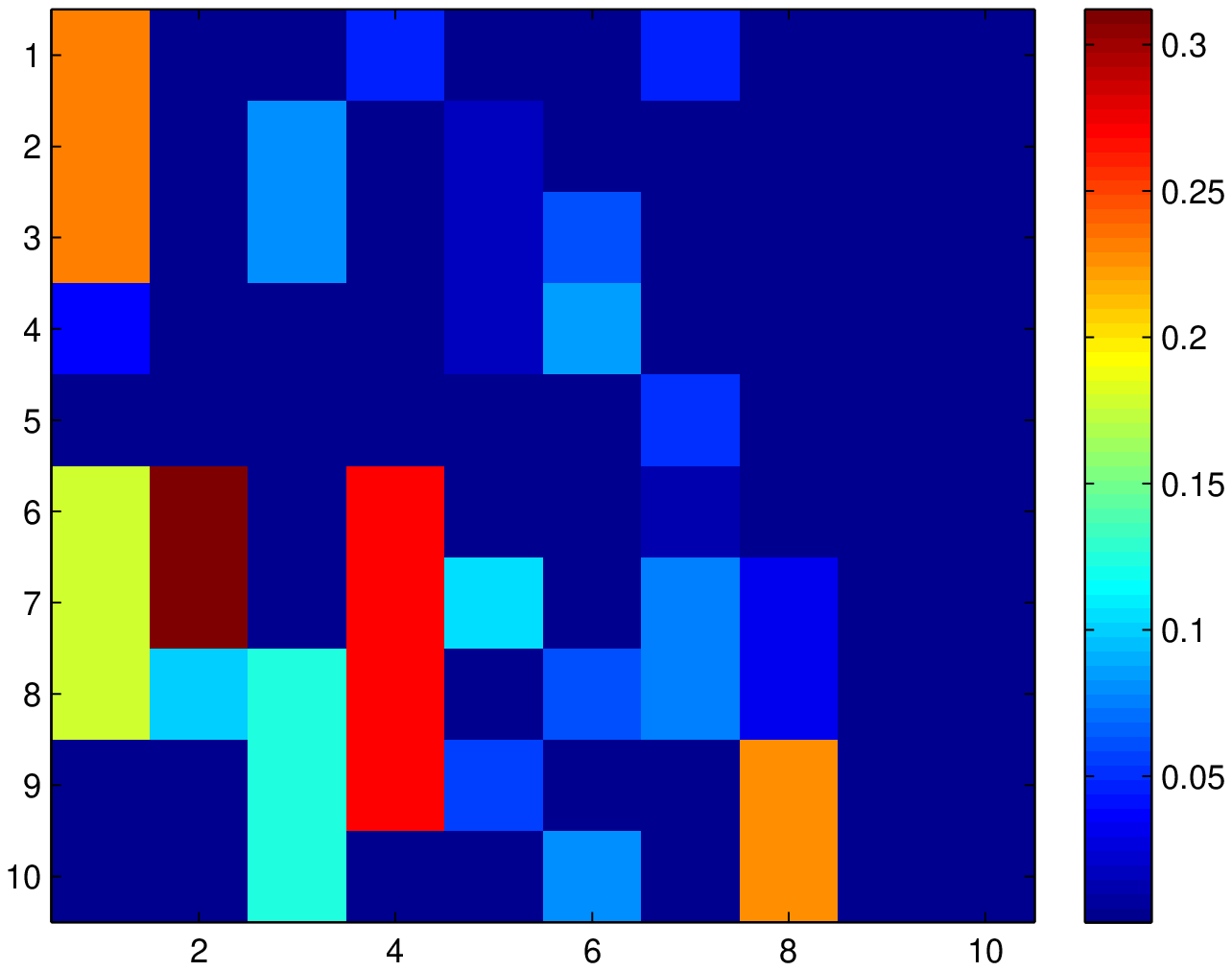}}\\ 
\subfloat[][]{\label{sparseheat-c}\includegraphics[width=2.5in]{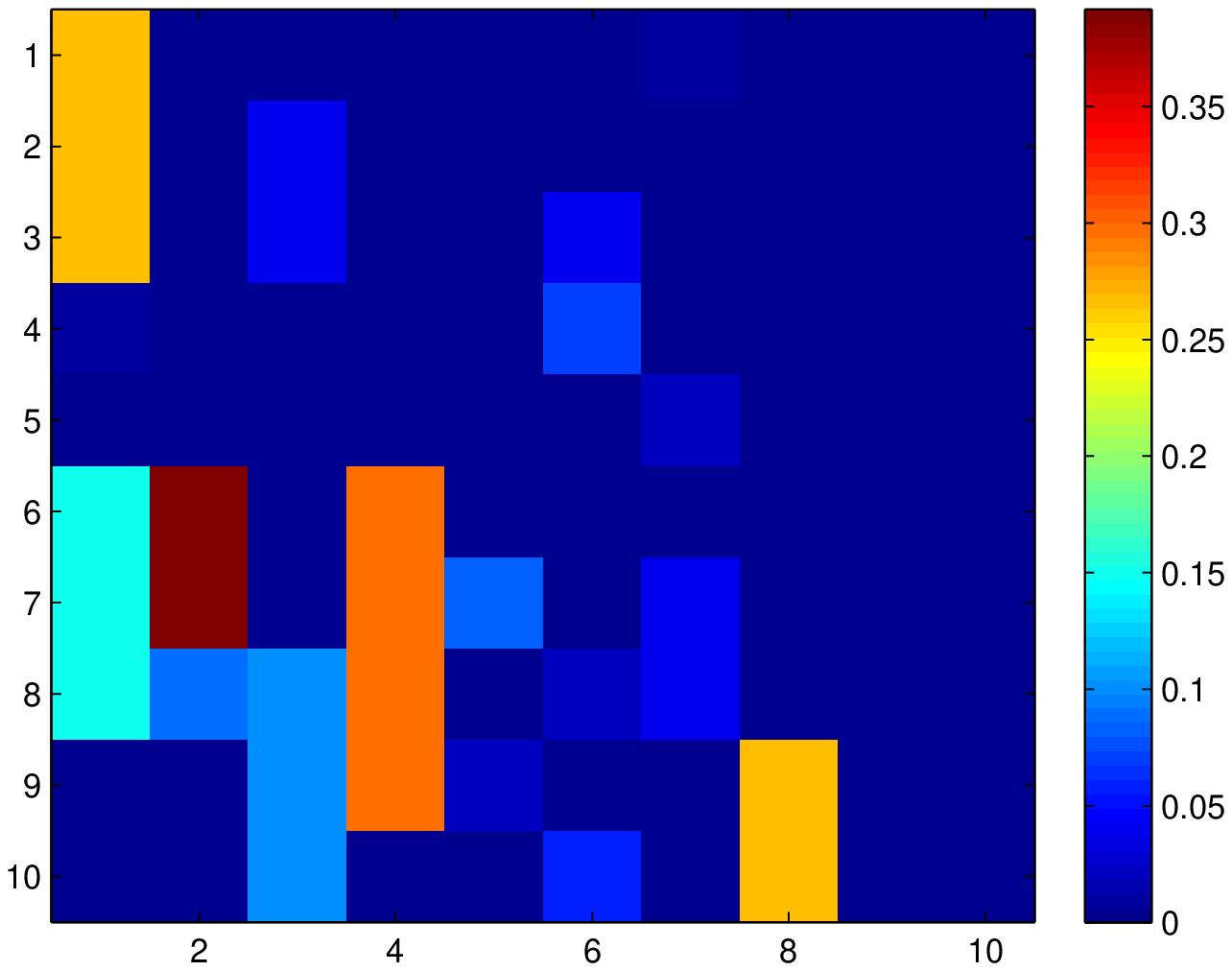}} 
\caption[Identifying responsible neural patches via $\ell_p$ norm.]{
\subref{sparseheat-a} Heat map using $\ell_{1/2}$ norm; 
\subref{sparseheat-b} Heat map using $\ell_{1/5}$ norm; 
\subref{sparseheat-c} Heat map using $\ell_{.0002}$ norm}
\label{sparseheat}%
\end{figure}

\begin{figure}%
\centering 
\subfloat[][]{\label{wristheat-a}\includegraphics[width=2.5in]{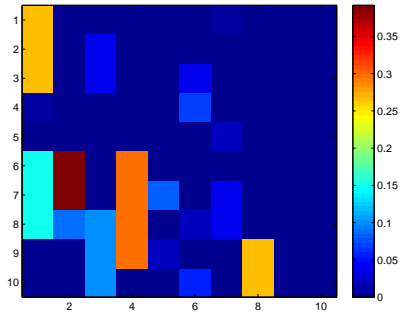}}\quad 
\subfloat[][]{\label{wristheat-b}\includegraphics[width=2.5in]{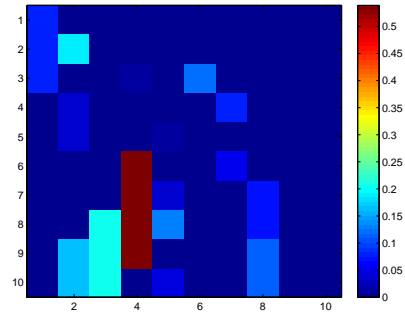}}\\ 
\subfloat[][]{\label{wristheat-c}\includegraphics[width=2.5in]{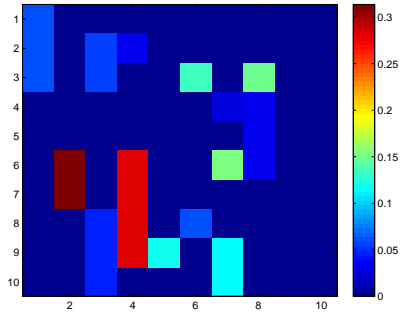}} 
\caption[Identifying responsible neural patches via $\ell_p$ norm.]{
\subref{wristheat-a} Heat map using $\ell_{.0002}$ norm for wrist up; 
\subref{wristheat-b} Heat map using $\ell_{.0002}$ norm for wrist down; 
\subref{wristheat-c} Heat map using $\ell_{.0002}$ norm for wrist right}
\label{wristheat}%
\end{figure}

\begin{figure}
\centering
\subfloat[][]{\label{diffheat-a}\includegraphics[width=2.5in]{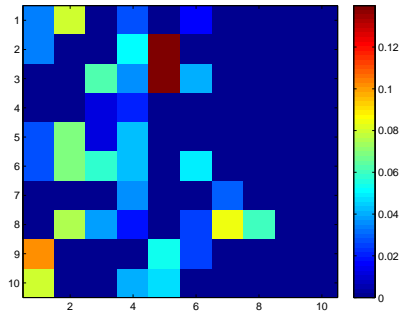}}\quad
\subfloat[][]{\label{diffheat-b}\includegraphics[width=2.5in]{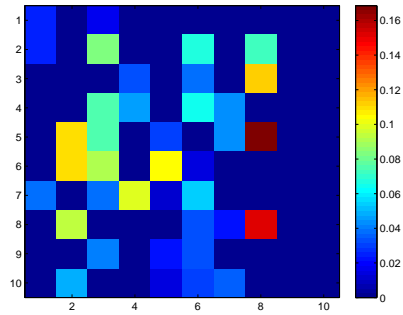} }
\caption[]{\subref{diffheat-a} Heat map for wrist up versus for $\alpha =1, p=1$;
\subref{diffheat-b} Heat map for shoulder up for $\alpha =1, p=1$.}
 \label{diffheat}
\end{figure}

Once the data has been classified, one can determine the responsible neurons for a given movement, based
on the classification.  To accomplish this, we form a heat map of the weights, $w$, corresponding to each neuron.
Each unit on the implanted patch consists of up to four neurons.  For this particular problem,
we are given a $10\times 10$ matrix, $M$, representing the physical location of the neurons with respect to the implanted patch.  Using
this matrix $M$ we map the weights, $w$, to the corresponding neuron's location on the patch.  A depiction of the heat map for several cases
can be seen in Figures \ref{sparseheat}, \ref{wristheat}, and \ref{diffheat}.
As can be seen in Figure \ref{sparseheat}, utilizing the $\ell_p$ norm for $p\ll 1$ allows one to identify
the most responsible neural regions for a specific movement.  As $p\to 0$ the number of identified neural regions is reduced.
For different wrist movements, the same region is identified, however when comparing shoulder versus wrist movements
different regions are identified with high activity, as can be seen in Figures \ref{wristheat} and \ref{diffheat}, respectively.

\clearpage
\section{Conclusion}
To summarize, we have provided a very generalized and easy to implement classification algorithm, based on the nonsmooth Tikhonov regularization
and the least squares formulation in \cite{psvm}.  The method is 
robust and works well even if the problem is severely ill-posed.  The algorithm not only provides the classifier, but is capable of detection
using the force as discussed in the results section.  By utilizing sparsity, one can minimize storage while determing the
responsible neurons for this particular application.
Future research includes more theoretical analysis of the methods proposed, as well as an extension to 
 nonlinear classification problems via reproducing kernel Hilbert space techniques.

\bibliography{neural}{}
\bibliographystyle{plain}
\nocite{psvm}
\nocite{newchoice}
\nocite{regparam}
\nocite{engl}
\nocite{cortes}
\nocite{bouton}
\clearpage

\appendix
\section{Convergence of the Iterative Algorithm}
\label{appenda}

Multiplying \eqref{iterate} by $w^{k+1}-w^k$, we obtain
$$\begin{array}{l}
\ds \frac{1}{2}(Hw^{k+1},w^{k+1})-(Hw^k,w^k)+(H(w^{k+1}-w^k),w^{k+1}-w^k))
\\ \\
\ds\quad +\frac{\beta p}{\max(\varepsilon^{2-p},|w^k|^{2-p})}\,\frac{1}{2}\,(|w^{k+1}|^2-|w^k|^2+
|w^{k+1}-w^k|^2)+(e,w^{k+1}-w^k).
\end{array} $$
Define the function
$$
\Psi_\varepsilon(x)= \left\{\begin{array}{ll} \ds \frac{p}{2}\frac{x}{\max\{\varepsilon^{2-p},|x|^{2-p}\}}+(1-\frac{p}{2})\,\varepsilon^p
& \quad x \le \varepsilon^2
\\ \\ w^{\frac{p}{2}}  &\quad x \ge \varepsilon^2 \end{array} \right.
$$
for $x\ge 0$ and
\begin{equation}  \label{eps}
J_\varepsilon(w)=\frac{1}{2}\,|Hw-e|^2+\Psi_\varepsilon(|w|^2).
\end{equation}
Then,
$$
\frac{1}{\max(\varepsilon^{2-p},|w^k|^{2-p})}\,\frac{p}{2}\,(|w^{k+1}|^2-|w^k|^2)=
\Psi^\prime_\varepsilon(|w^k|^2)(|w^{k+1}|^2-|w^k|^2)
$$
Since $x \to \Psi_\varepsilon(x)$ is concave, we have
$$
\Psi_\varepsilon(|w^{k+1}|^2)-\Psi_\varepsilon(|w^k|^2)
-\frac{1}{\max(\varepsilon^{2-p},|w^k|^{2-p})}\,\frac{p}{2}\,(|w^{k+1}|^2-|w^k|^2)\le 0
$$
and thus
\begin{equation} \label{est}
J_\varepsilon(w^{k+1})+\frac{1}{2}(H(w^{k+1}-w^k),w^{k+1}-w^k))
+\frac{\beta p}{\max(\varepsilon^{2-p},|w^k|^{2-p})}\,\frac{1}{2}\,|w^{k+1}-w^k|^2\le J_\varepsilon(w^k)
\end{equation}
shows that $J_\varepsilon$ is nonincreasing.  Now, we give the following result.

\begin{theorem} For $\varepsilon>0$ let $\{w_k\}$ be generated by \eqref{iterate}. Then,
$J_\varepsilon(w^k)$ is monotonically non increasing and $w_k$ converges to the minimizer of $J_\varepsilon$
defined by \eqref{eps}. \vspace{2mm}
\end{theorem}

\begin{proof} The monotonicity of $J_\varepsilon$ has already been shown.  Thus, we show that $\{w^k\}$ converges to the minimizer of $J_\varepsilon$.
It follows from \eqref{est} that $|w^k|_\infty <\infty$
and
$$
\sum^\infty_{k=0}|w^{k+1}-w^k|_2^2 < \infty
$$
and thus there exists subsequence of $\{w^k\}$ and $w^*\in\ell^p$ such that
$$
\lim_{k \to \infty} w_k=\lim_{k\to\infty} w^{k+1}=w^*.
$$
It follows from \eqref{iterate}
$$
H^*Hw^*+\frac{\beta p}{\max(\varepsilon^{2-p},|w^*|^{2-p})}\,w^*=H^*e,
$$
i.e., $w^*$ minimizes $J_\varepsilon$.
\end{proof}

\end{document}